\documentclass[a4paper,12pt,reqno]{amsart}
\usepackage{amssymb}
\usepackage{latexsym}
\usepackage{amsmath}
\usepackage[english]{babel}
\usepackage{amsfonts}
\usepackage{enumerate}
\usepackage[all]{xy}
\usepackage[dvipdfm]{pict2e}
\usepackage[dvips]{graphicx, color}\usepackage{verbatim}
\CompileMatrices

\newlength{\widthuparrow}
\newcommand{\ind}{\operatorname{ind}}

\newcommand{\lbr}{\begin{bmatrix}}
\newcommand{\rbr}{\end{bmatrix}}
\newcommand{\Ker}{\operatorname{Ker}}
\newcommand{\I}{\operatorname{Im}}
\newcommand{\im}{\operatorname{im}}
\newcommand{\Hom}{\operatorname{Hom}}
\newcommand{\Ext}{\operatorname{Ext}}

\newcommand{\rd}{\operatorname{d}}
\newcommand{\St}{\operatorname{St}}
\newcommand{\soc}{\operatorname{soc}}
\newcommand{\rad}{\operatorname{rad}}
\newcommand{\hd}{\operatorname{hd}}

\newcommand{\ch}{\operatorname{ch}}
\newcommand{\pr}{\operatorname{pr}}

\usepackage{bbm} 
\newcommand\cA{\mathcal A}
\newcommand\cC{\mathcal C}
\newcommand\C{\mathbb C}
\newcommand\Z{\mathbb Z}

\newcommand\N{\mathbb N}
\newcommand\Q{\mathbb Q}

\newcommand\cc{\mathcal C}

\usepackage{amsthm}
\newtheorem{thm}{Theorem}[section]

\newtheorem{prop}[thm]{Proposition}
\newtheorem{cor}[thm]{Corollary}
\newtheorem{lemma}[thm]{Lemma}

\theoremstyle{definition}
\newtheorem{defn}[thm]{Definition}
\newtheorem{rem}[thm]{Remark}
\newtheorem{conj}[thm]{Conjecture}

\numberwithin{equation}{section}

\title{Rigidity of tilting modules}
\author{Henning Haahr Andersen and Masaharu Kaneda }

\address{HHA: Department of Mathematics, University of Aarhus, Building 530, Ny Munkegade,
8000  Aarhus C, Denmark}
\email{mathha@imf.au.dk}
\thanks
{MK supported in part
by JSPS Grant in Aid
for Scientific Research}

\address{MK: Department of Mathematics, Osaka City University,
Osaka 558-8585 Japan}
\email{kaneda@sci.osaka-cu.ac.jp}

\begin{document}

\maketitle

\section*{Abstract}
Let $U_q$ denote the quantum group associated with a finite dimensional semisimple Lie algebra.
Assume that $q$ is a complex root of unity of odd order and that $U_q$ is 
obtained via
Lusztig's $q$-divided powers construction. We prove that all regular projective (tilting) modules for $U_q$ are 
rigid, i.e., have identical radical and socle filtrations. Moreover, we obtain the same 
for a large class of Weyl modules for $U_q$. On the other hand, we give
examples of non-rigid indecomposable tilting modules as well as non-rigid Weyl modules. These examples are for type $B_2$ 
and in this case as well as for type $A_2$ we calculate explicitly the Loewy structure for all regular Weyl modules.

We also demonstrate that these results carry over to the modular case when the highest weights in question are
in the so-called Jantzen region. At the same time we show by examples that as soon as we leave this region 
non-rigid tilting modules do occur.

\noindent {17B37, 20G05}

\section{Introduction}

This paper concerns the structure of tilting modules for quantum groups at 
complex roots of unity together with the corresponding modules for semisimple
algebraic groups over fields of positive characteristics.
Our aim is to determine the Loewy structure of these modules. In the quantum
case we prove that all indecomposable projective tilting modules with regular
highest weights are rigid (i.e. have identical radical and socle series). We also give examples (for type $B_2$) showing
that rigidity fails for some non-projective indecomposable tilting modules.

Our proof makes use of the Lusztig conjecture for the irreducible characters of
quantum groups at complex roots of $1$. This conjecture is a theorem, see \cite{KL}, 
\cite {L94} and \cite {KT}. More recently, alternative (and shorter) proofs have
been given in  \cite{ABG} and \cite{F1}. This result implies that also the
characters of all indecomposable projective tilting modules are known. Moreover,
due to \cite{Soe98} we know the characters of all indecomposable tilting modules
and we take advantage of this when working out our $B_2$-examples.

In the modular case the Lusztig conjecture on the irreducible characters for a
semisimple algebraic group is only proved for large primes \cite{AJS} (see also
the recent work by Fiebig \cite{F2} giving an explicit bound for how large the
prime suffices to be). Still worse: the characters of indecomposable tilting modules
are completely unknown - there is not even a conjecture. In this case we proved
many years ago \cite{AK} that if we replace the semisimple algebraic group $G$
by $G_1T$ where $G_1$ is the Frobenius kernel and $T$ is a maximal torus then
the indecomposable $G_1T$-modules corresponding to regular highest weights are
rigid and we determined the composition factor multiplicities of their Loewy
layers. To establish this result we assumed the Lusztig conjecture on
irreducible characters. With the same assumption and also assuming $p$ to be
bigger than $3(h-1)$ with $h$ denoting the Coxeter number for $G$ it is easy to
deduce that indecomposable tilting modules with regular highest weight above the
Steinberg weight and below the upper bound for the Jantzen region (see Theorem
4.7 below for the exact statement) are also rigid. Our $B_2$-example mentioned
above also works in this modular case. It shows that we cannot expect to relax
the condition that the weight should lie above the Steinberg weight. We give
another example - this time for type $A_2$ - of a non-rigid indecomposable
tilting module whose highest weight lies slightly above the upper bound in this
theorem. This illustrates the difference between the behavior of tilting modules
in the modular and the quantum case. At the same time the example proves the
suspicion expressed in \cite{DM} that certain $SL_3$-modules are non-rigid (not
only in very small characteristics as in loc.cit. but for all odd primes).

Tilting modules are characterized by having filtrations by both Weyl modules and
dual Weyl modules. These modules play naturally a key role in our treatment. Along 
the way we also establish their rigidity for certain weights (see
Theorem 3.24 and Remark 3.25). Although we believe these
conditions to be too restrictive our $B_2$-examples in Section 5.9 show that
some restrictions are definitely needed. In type $A_2$, however, we
find that all regular quantum Weyl modules are rigid, see Corollary 5.5.

Among other consequences of our work we mention that in the quantum case we get
a universal upper bound on the Loewy length of arbitrary finite dimensional
modules. In the modular case the same bound holds but only for modules with
weights below the upper bound mentioned above. 

The paper is organized as follows: We introduce the basic definitions and recall
some standard facts about rigidity and tilting modules in Section 2. Then we
treat the quantum case in Section 3, the modular case in Section 4, give our 
$B_2$-examples in Section 5, and finish with the mentioned $A_2$-computations in
Section 6. In a short appendix we prove a result which we need in Section 5 to determine 
the socles of certain Weyl modules.
   
\section{Notations and some basic properties}

\subsection{Roots and Weights}
Let $R$ denote a (finite) root system in the Euclidian space $V$ and pick a set
of positive roots $R^+$. We denote by $X =
\{\lambda \in V \mid \langle \lambda, \alpha^\vee \rangle \in \Z \text { for all
} \alpha \in R\}$ the set of
integral weights and by $X^+ =
\{\lambda \in X \mid \langle \lambda, \alpha^\vee \rangle \geq \text { for all }
\alpha \in R^+ \}$ the set of dominant weights. Here $\alpha^\vee$ is the
coroot associated to $\alpha$.

The Weyl group for $R$ is generated by $s_\alpha$ with $\alpha$ running through
$R$. As usual we shall shift the usual action of $W$ on $V$ and $X$
by $-\rho$ where $\rho$ is half the sum of the positive roots. We write $w \cdot
\lambda = w(\lambda + \rho) - \rho,\; w \in W, \: \lambda \in X$. 

We let $A = (a_{ij})_{i,j = 1, \cdots , n}$ denote the Cartan matrix for $R$.
We pick $(d_i)_{i = 1, \cdots , n}$ minimal such
that $(d_i a_{ij})_{i,j = 1, \cdots , n}$ is symmetric.

\subsection{The modular case}
      
Let $k$ denote an algebraically closed field of characteristic $p > 0$ and
let $G$ be a connected reductive algebraic group over $k$. Fix a maximal torus
$T$ in $G$ and assume that $R$ is the root system associated to $(G, T)$.
To identify $X$ with the character group of $T$,
we will assume $G$ to be simply connected.

In the category $\mathcal C$ of finite dimensional $G$-modules we have four important
modules associated with a dominant weight $\lambda$. First we have the
irreducible module $L(\lambda)$ with highest weight $\lambda$. This is the
unique simple quotient of the Weyl module $\Delta (\lambda)$ as well as the
unique simple submodule of the dual Weyl module $\nabla (\lambda)$. Finally, we
have the indecomposable tilting module $T(\lambda)$. This module has a Weyl
filtration starting with $\Delta (\lambda)$ and a dual Weyl filtration ending
with $\nabla (\lambda)$. All other quotients in these two filtrations of $T
(\lambda)$ have weights strictly less than $\lambda$. So $\lambda$ is the unique
highest weight in all four modules $L(\lambda), \Delta (\lambda),
\nabla(\lambda) $, and $T(\lambda)$, and it occurs with multiplicity $1$ in each of
them.

\subsection{The quantum case}

Let $v$ denote an indeterminate and set $U_v$ equal to the quantum group over 
$\Q(v)$ associated to the Cartan matrix $A$ of $R$. The generators
of this $\Q(v)$-algebra are $E_i, F_i,$ and $K_i^{\pm 1}, i=1, 2, \cdots, n$. 
We set $U_v^+$, respectively $U_v^-$, respectively $U_v^0$ equal to the 
subalgebra
generated by the $E_i$'s, respectively $F_i$'s, respectively $K_i^{\pm 1}$'s. 
Then $U_v = U_v^- U_v^0 U_v^+$.

Set $\mathcal A = \Z[v,v^{-1}]$ and consider the Lusztig $\mathcal A$-form of 
$U_v$. This is the $\mathcal A$-subalgebra of $U_v$ generated by 
$E_i^{(r)}, F_i^{(r)},
K_i^{\pm 1}$, $i= 1, 2, \cdots , n, \; r \in \N$. Here the divided powers 
$E_i^{(r)}$ and $F_i^{(r)}$ are defined as follows: if $m \in \Z$ and $d \in \N$ then the 
Gaussian number $[m]_d$ is given by $[m]_d = \frac{v^{dm} - v^{-dm}}{v^d - v^{-d}}$ and for 
$m>0$ we set $[m]_d! = [m]_d [m-1]_d \cdots [1]_d$. Then $E_i^{(r)} =
\frac {E_i^r}{[r]_{d_i}!}$. We define $F_i^{(r)}$ by the same recipe.

Again we have subalgebras  $U_{\mathcal A}^+, U_{\mathcal A}^-, U_{\mathcal A}^0$ such that 
$U_{\mathcal A} = U_{\mathcal A}^- U_{\mathcal A}^0 U_{\mathcal A}^+$. 
Here $U_{\mathcal A}^+ $, respectively $U_{\mathcal A}^-$, is the $\mathcal A$-subalgebra of $U_{\mathcal A}$ generated by the $E_i^{(r)}$'s, respectively
$F_i^{(r)}$'s whereas $U_{\mathcal A}^0$ is generated by the $K_i^{\pm 1}$ together with the following elements 
$$\left[\begin{smallmatrix} K_i \\t 
\end{smallmatrix}\right]
=\prod^t_{j={1}} \frac{K_i v^{d_i(1-j)}
-K^{-1}_i v^{d_i(j-1)}}{v^{d_ij} - v^{- d_i
j}};\; t \in \N.$$

Let now $q \in \C$ be a root of unity of order $l$. We assume that $l$ is prime 
to all non-zero entries in $A$. In particular $l$ is odd. 
Then we make $\C$ into an $\mathcal A$-algebra by specializing  $v$ to $q$. The quantum group we want to consider is then 
$$ U_q = U_{\mathcal A} \otimes \C.$$
We set $U_q^+ = U_{\mathcal A}^+ \otimes \C, U_q^- = 
U_{\mathcal A}^- \otimes \C$ and $U_q^0 = U_{\mathcal A}^0 \otimes \C$.  
Again we have $U_q = U_q^- U_q^0 U_q^+$.

Let $\lambda \in X$. Then we have a character $\chi_\lambda: U_q^0 \rightarrow \C$ given by $\chi_\lambda(K_i) = 
q^{d_i\lambda_i}$ and $\chi_\lambda(\left[\begin{smallmatrix} K_i \\t 
\end{smallmatrix}\right]) = \left[\begin{smallmatrix} \lambda_i \\t 
\end{smallmatrix}\right]_{d_i}$ , see \cite{L90}. Here $\lambda_i = \langle \lambda, \alpha^\vee_i \rangle$ is the $i$'th coordinate of $\lambda$ (with $\alpha_1, \alpha_2, \cdots ,
\alpha_n$ denoting the simple roots) and for general $m\in \Z, t, d \in \N$ we write $\left[\begin{smallmatrix} m \\t 
\end{smallmatrix}\right]_{d} \in \C$ for the Gaussian binomial coefficient  $\left[\begin{smallmatrix} m \\t 
\end{smallmatrix}\right]_{d} = \frac{[m]_d!}{[t]_d! [m-t]_d!} \in \mathcal A$ evaluated at $q$.

If $M$ is a $U_q^0$-module then the weight space $M_\lambda$ is defined by $M_\lambda = \{m\in M \mid u m = \chi_\lambda (u) m \text { for all } u \in U_q^0\}$. We
say that $M$ has type $\bf 1$ if $M = \oplus_{\lambda \in X} M_\lambda$. 

In the following all modules we consider will be finite dimensional of type $\bf 1$. The category of such $U_q$-modules will be denoted $\mathcal C_q$. (In several
previous papers this notation was used for the bigger category of integral $U_q$-modules. However, we will only need finite dimensional modules for our purposes
here).

In $\mathcal C_q$ we have - in complete analogy with the modular situation described above - the following four modules attached to a weight $\lambda \in X^+$: 
The simple module $L_q(\lambda)$, the Weyl module $\Delta_q (\lambda)$, the dual Weyl module $\nabla_q(\lambda)$ and the indecomposable tilting module $T_q(\lambda)$.

Each of these modules have highest weight $\lambda$. The dual Weyl module $\nabla_q(\lambda)$ may be obtained via induction from the $1$-dimensional
$U_q^-U_q^0$-module given by the character $\chi_\lambda$ (extended trivially to $U_q^-U_q^0$). 
Then we can define $\Delta_q(\lambda)$ as the dual of $\nabla_q(\lambda)$. 
For the dual operation
we will employ an involutive antiautomorphism 
$\omega S$
of
$U_q$,
where $S$ is the antipode of
$U_q$
and $\omega$ is an involutive automorphism of
$U_q$
exchanging each $E_i$ and $F_i$ while sending
$K_i$ to $K_i^{-1}$.
For
$M\in\mathcal C_q$
we let
$M^{\omega S}$ denote the $\mathbb{C}$-linear
dual of
$M$ with $U_q$ acting through $\omega S$.
So we set  $\Delta_q(\lambda)=\nabla_q(\lambda)^{\omega S}$ and $L_q(\lambda)$
 is realized as the socle of $\nabla_q(\lambda)$ (or equivalently the head of $\Delta_q(\lambda)$). The tilting module $T_q(\lambda)$ contains $\Delta_q(\lambda)$
 and surjects onto $\nabla_q(\lambda)$.

\subsection{Rigid modules}
Let $M$ denote a module in either $\mathcal C$ or $\mathcal C_q$. The radical
series of $M$ is then the series of submodules
$$ 0 = \rad^r M \subset \rad^{r-1} M \subset \cdots \subset \rad^1 M \subset M$$
where by the radical $\rad ^1 M$ of $M$ we mean the smallest submodule whose
quotient is semisimple. We also often write just $\rad M$ instead of $\rad^1 M$
and we set $\rad^{i+1} M = \rad (\rad^i M)$. The quotient $M/\rad M$ is called the 
head of $M$ and denoted $\hd M$. We set $\rad_i M = \rad^i M/\rad^{i+1} M$ and 
call this the $i$-th radical layer of $M$.

Similarly, the socle series of $M$ is the series of submodules 
$$ 0 = \soc^0 M \subset \soc^1 M \subset ^2 M \subset \cdots \subset \soc^r M =
M.$$
Here $\soc^1 M$ is the largest semisimple submodule of $M$, often denoted just
$\soc M$, and $\soc^{i+1} M = \pi^{-1}(\soc (M/\soc^ i M))$ with $\pi$ being
the projection $M \rightarrow M/\soc^i M$. The $i$-th socle layer of $M$ is 
$\soc_i M = \soc^{i} M/\soc ^{i-1} M$. 

As indicated the length $r$ of the radical series for $M$ coincides with the
length of the socle series. This common number is called the Loewy length and
denoted $ll M$.

\begin{defn}
The module $M$ is rigid if its radical series coincides with it socle series,
i.e., if $\rad^i M = \soc^{ll M - i } M$.
\end{defn}

A Loewy series for $M$ is a filtration $ 0 = M^r \subset M^{r-1} \subset  \cdots
\subset M^1 \subset M^0 = M$ with the property that $r = ll M$ and all quotients
$M^i/M^{i+1} M$ are semisimple. By definition we have $\rad ^i M \subset M^i$
and
$M^i \subset \soc^{ll M - i } M$. So $M$ is rigid if and only if it has a unique Loewy
series.

\begin{remark}
Clearly, the dual of the radical series for $M$ gives us the socle series for the dual module $M^*$. 
Explicitly, we have $(\rad ^i M)^* \simeq M^*/\soc^i (M^*)$ for all $i$.
\end{remark}

\section{Loewy structure of modules for quantum groups}

In this section we study the Loewy structure of some of the tilting modules in 
the category $\mathcal C_q$ from Section 2. For simplicity we assume that our
root system $R$ is irreducible. Recall that $q$ is a complex root of unity of 
odd order $l$ (with $l$ also being prime to $3$ if $R$ is of type $G_2$).

 \subsection{Projective modules}
 Consider the weight $(l-1)\rho \in X^+$. In this case the strong linkage principle \cite[7.6]{APW91} gives
 $$ L_q((l-1)\rho) = \nabla_q((l-1)\rho) = \Delta_q((l-1)\rho) = T_q((l-1)\rho).$$
We denote this special module $\St_q$ and call it the Steinberg module. It has (again by the strong linkage principle) the property
 $$ \St_q \text { is injective (and projective) in  } \mathcal C_q.$$
 
 Note that $\St_q$ is selfdual. Hence the trivial module $k$ is a submodule of $\St_q \otimes \St_q$ and any module $M \in \cc_q$ is embedded into $M \otimes \St_q
 \otimes \St_q$, i.e., we have
 $$ \cc_q \text { has enough injectives and enough projectives}.$$
 
 We set $X_l = \{\nu \in X^+ \mid \langle \nu, \alpha^\vee \rangle < l \text { for all simple roots } \alpha \}$ 
 and write for each $\lambda \in X$, $\lambda = \lambda^0 + l \lambda^1$ with $\lambda^0 \in X_l$. Then we set $\tilde \lambda = 2(l-1)\rho +w_0\lambda^0 + l \lambda^1$,
 where $w_0\in W$ is the longest element
 such that
 $w_0R^+=-R^+$.

 \begin{prop} 
\cite[5.8]{A92}
 Let $\lambda \in X^+$. Then 
  $T_q(\tilde \lambda)$ is the injective envelope (and the projective cover) of $L_q(\lambda)$.

  \end{prop}
  

 \begin{cor}
 Let $\lambda \in X^+$. Then $L_q(\lambda) = \soc \Delta_q(\tilde \lambda) = \hd \nabla_q(\tilde \lambda)$.
 \end{cor}
 
 \begin{pf}
 Note that $\Delta_q(\tilde \lambda)$ is a submodule of $T_q(\tilde \lambda)$ and $\nabla_q(\tilde \lambda)$ is a quotient. Now apply Proposition 3.1.
 \end{pf}
 
 \begin{rem}
 This corollary proves that if $\nu \in (l-1)\rho + X^+$ then $\nabla_q(\nu)$ has simple head. This is not true for all
 $\nu \in X^+$ as our $B_2$-examples in Section 5 will illustrate. In fact, in that case we determine the heads of all dual Weyl modules, see Remark 7.3(2).
 \end{rem}
 
 \subsection{Alcoves and wall crossings}
 The bottom alcove $C$ in $X^+$ is defined by
 $$ C = \{\lambda \in X^+ \mid \langle \lambda + \rho, \alpha^\vee \rangle < l \text { for all } \alpha \in R^+ \}$$
 and its closure $\bar C$ is 
 $$\bar C = \{\lambda \in X \mid 0 \leq \langle \lambda + \rho, \alpha^\vee \rangle \leq l \text { for all } \alpha \in R^+ \}.$$
 
 The affine Weyl group $W_l$ may then be defined as the group generated by the reflections in the walls of $C$. An arbitrary alcove is $w \cdot C$ with $w \in W_l$. 
 
 A weight $\lambda \in X$ is called $l$-regular provided $\langle \lambda + \rho, \alpha^\vee \rangle \not \equiv 0 \; \mod l  \text { for all } \alpha \in R^+$.
 Equivalently, $\lambda$ is $l$-regular iff $\lambda$ belongs to some alcove in $X$.
 
 We now let $S_l$ denote the set of reflections in the walls of $C$. Then $S_l$ is a set of generators for $W_l$. If $s \in S_l$ then any alcove $A$ has a unique wall 
 which is in the $W_l$-orbit of this wall. We call this the
 $s$-wall of $A$ and we write $As$ for the alcove obtained by reflecting $A$ in its $s$-wall.

Denote by $\mathcal A^+$ the set of alcoves in $X^+$. For later use we also introduce $\mathcal A^{++}$ as the set of alcoves in $l\rho + X^+$. 
 
 To each $\lambda \in \bar C$ we associate the block $\mathcal B_\lambda \subset \mathcal C_q$ consisting of those $M \in \cc_q$ whose compositions factors all have
 highest weights belonging to the orbit $W_l \cdot \lambda$. Then 
 $$ \cc _q = \oplus _{\lambda \in \bar C} \mathcal B_\lambda.$$
 We denote the corresponding projection $\cc _q \rightarrow \mathcal B_\lambda$ by $\pr_\lambda$.
 
 Let now also $\mu \in \bar C$. Then we have a translation functor $T_\mu^\lambda : \mathcal B_\mu \rightarrow \mathcal B_\lambda$ defined by
\[
T_\mu^\lambda M = \pr_\lambda (M \otimes T_q(x(\lambda - \mu))).
\]
 Here $x \in W$ is chosen such that $x(\lambda - \mu) \in X^+$.
 We have for any $V\in\mathcal C_q$
 \[
(M\otimes V)^{\omega S}\simeq
M^{\omega S}\otimes V^{\omega S},
 \]  
 and hence 
 \[
(T_\lambda^\mu M)^{\omega S}\simeq
T_\lambda^\mu(M^{\omega S}).
 \]
 
 \begin{remark} Usually the translation functors $T_\mu^\lambda$ are defined by tensoring with simple modules instead of tilting modules. However, the outcome only
 depends on the extremal weights (those in the $W$-orbit of the highest weight) of the module used. Our choice makes it clear that translation functors take tilting
 modules to tilting modules, see \cite{P},
 \cite{X},
 \cite{K98}.
 \end{remark}
 
 \subsection{Loewy lengths}

 In this and the following subsections we fix an $l$-regular weight $\lambda \in C$. If $A$ is any other alcove then we denote by $\Delta_q(A)$ the Weyl module with highest 
 weight in $W_l\cdot \lambda \cap A$. Similarly, we define $L_q(A), \nabla_q(A)$ and $T_q(A)$. 
 
 Let $s \in S_l$ and choose a weight $\mu$ in the interior of the $s$-wall of $C$. Then the composite $T_\mu^\lambda \circ T_\lambda^\mu : \mathcal B_\lambda
 \rightarrow \mathcal B_\lambda$ is denoted $\theta_s$. It is sometimes called the wall-crossing functor.
 
 If $\nu \in \bar C$ and $A$ is an arbitrary alcove then we let $\nu_A$ denote the weight given by $\{\nu_A\} = W_l\cdot \nu \cap \bar A$. If $A$ is determined by
 $A = \{\eta \in X \mid n_\alpha l < \langle \eta + \rho, \alpha^\vee \rangle < (n_\alpha + 1) l \text { for all } \alpha \in R^+ \}$ where $n_\alpha \in \Z$ 
 then the upper
 closure $\hat A$ of $A$ is defined by 
 $\hat A =  \{\eta \in X \mid n_\alpha l < \langle \eta + \rho, \alpha^\vee \rangle \leq (n_\alpha + 1) l \text { for all } \alpha \in R^+ \}$.
 
 The following proposition is the analogue of well-known modular results, see Proposition 7.15 in \cite{RAG}
 \begin{prop}
 Let $A \in \mathcal A^+$. 
 \begin{enumerate}
 \item[i)] For any $\nu \in \bar C$ we have $T_\lambda^\nu \nabla_q(A) = \nabla_q (\nu_A)$ and \\ $T_\lambda^\nu L_q(A) =
 \begin{cases} L_q(\nu_A) & \text { if } \nu_A \in \hat A,\\ 0 & \text { otherwise.} \end{cases}$
 \end{enumerate}
Assume now also that $As \subset X^+$. Then
\begin{enumerate} 
 \item[ii)] $\theta_s \nabla_q(A) \simeq \theta_s \nabla_q(As)$.
 \item [iii)] Suppose $As < A$. Then we have a non-split exact sequence 
 $0 \rightarrow \nabla_q(As) \rightarrow \theta_s \nabla_q(A) \rightarrow \nabla_q(A) \rightarrow 0.$
 \end{enumerate}
 If $As > A$ then there is a corresponding sequence with the roles of $A$ and $As$ reversed. 

 \end{prop}
 
 \begin{cor}
 With notation as in the above proposition and assuming $As < A$ we have $\Hom_{U_q}(\nabla_q(A), \nabla_q(As)) \simeq \C$.
 \end{cor}
 
 \begin{pf} 
 The functors $T_\lambda^\mu$ and $T_\mu^\lambda$ are adjoint. Hence we get 
 \begin{align*}
 \Hom_{U_q}(\nabla_q(A), \theta_s \nabla_q(A)) &\simeq \Hom_{U_q}(T_\lambda^\mu
 \nabla_q(As), T_\lambda^\mu \nabla_q(A))
 \\
 &\simeq \Hom_{U_q}(
 \nabla_q(\mu_A), \nabla_q(\mu_A))
 \simeq \C.
 \end{align*}
The non-splitness of the exact sequence in Proposition 3.4 iii) implies that any non-zero element in $\Hom_{U_q}(\nabla_q(A),
\theta_s\nabla_q(A) )$ is in fact
 in $\Hom_{U_q} (\nabla_q(A),
 \nabla_q(As))$. The corollary follows.
 \end{pf}
 
 \begin{rem}
 The proof of this corollary shows that if $\phi \in$
 \linebreak 
 $\Hom_{U_q}(\nabla_q(A), \nabla_q(As))$ is non-zero then $T_\lambda^\mu \phi$ is an isomorphism.
 \end{rem}
 
 Set now $A^+ = C + l\rho$ and $A^- = w_0 \cdot C + l\rho$ and choose a sequence 
 $$ A_0 = A^+ > A_1 > A_2 > \cdots > A_N = A^-$$
 such that $(l-1)\rho \in \bar A_i$ for all $i$ and such that $A_i$ and $A_{i+1}$ share a common wall, say an $s_i$-wall, for which $A_is_i  = A_{i+1}$. The length of
 such a sequence is $N = \# R^+$. We choose $\mu_i$ in the interior of the $s_i$ wall of $C$.
 
 By Corollary 3.5 we have up to scalars a unique non-zero homomorphism $\varphi_{i+1} \in \Hom_{U_q}(\nabla_q(A_i), \nabla_q(A_{i+1}))$.
 
 \begin{prop}
 In the above notation we have that the composite $\varphi_N \circ \varphi_{N-1} \circ \cdots \circ \varphi_1 : \nabla_q(A^+) \rightarrow \nabla_q(A^-)$ is non-zero.
 \end{prop}
 
 \begin{pf} By Corollary 3.2 we have $\hd \nabla_q(A^+) = L_q(A^-)$. Noting that for each $i$ the $s_i$-wall of $A^-$ is in the upper closure of $A^-$ we get from Proposition 3.4
 i) that $T_\lambda^{\mu_i}L_q(A^-)$ is non-zero. If $L_q(A^-)$ were a composition factor of $\Ker (\varphi_1)$ then the kernel of $T_\lambda^{\mu_1} \varphi_1$ would
 be non-zero. But Remark 3.6 says that this homomorphism is an isomorphism. Similarly $L_q(A^-)$ cannot be a composition factor of the kernels of any of the other
 $\varphi_i$'s either. The proposition follows.
 \end{pf}
 
\begin{rem}
 It is easy to see that in fact $\Hom_{U_q}(\nabla_q(A^+), \nabla_q(A^-)) \simeq \C$. This means that any element in this space is a scalar multiple of the composite
 $\nabla_q(A^+) \rightarrow \hd \nabla_q(A^+) = L_q(A^-) = \soc \nabla_q(A^-)  \subset \nabla_q(A^-)$. Hence we have in fact $\I (\varphi_N \circ \varphi_{N-1} \circ
 \cdots \circ \varphi_1) = L_q(A^-)$.
 \end{rem}
 
 \begin{cor}
 $ll (\nabla_q(A^+)) \geq N+1$.
 \end{cor}
 
 \begin{pf} Since $L_q(A_i)$ is not a composition factor of $\nabla_q(A_{i+1})$ we see that $\varphi_{i+1}$ vanish on $\soc \nabla_q(A_i) = L_q(A_i)$. This gives
 the inequalities (the first is actually an equality, see Remark 3.8)
\begin{multline*}
 1 \leq ll(\I(\varphi_N \circ \varphi_{N-1} \circ\cdots \circ \varphi_1)) < ll(\I(\varphi_{N-1} \circ \varphi_{N-2} \circ\cdots \circ \varphi_1)) < \cdots
 \\ 
<
 ll(\I(\varphi_1)) < ll (\nabla_q(A^+))
 \end{multline*}
 and we are done.
 \end{pf}
 
 \begin{prop} For any alcove $A \in \mathcal A^+$ we have
\[
ll (T_q(A)) \geq 2 ll(\nabla_q(A)) - 1.
\]
 \end{prop}
 
 \begin{pf}
 Set $m = ll (\nabla_q(A))$. Taking the 
 $\omega S$-dual, we have also $m = ll (\Delta_q(A))$.
 Recall that $\nabla_q(A)$ is a quotient of $T_q(A)$ while $\Delta_q(A)$  is a submodule. The simple module $L_q(A)$ occurs with multiplicity $1$ in $T_q(A)$.
 Suppose $\rad^i T_q(A)$ does not contain $L_q(A)$ as a composition factor. Then the surjection from $T_q(A)$ to $\nabla_q(A)$ must factor through 
 $T_q(A)/\rad^iT_q(A)$ and hence $i = ll(T_q(A)/\rad^iT_q(A)) \geq ll (\nabla_q(A))= m$. We conclude that $L_q(A)$ must occur as composition factor of $\rad^{m-1}
 T_q(A)$. It follows that the submodule $\Delta_q(A)$ of $T_q(A)$ must in fact be a submodule of $\rad^{m-1} T_q(A)$. Therefore we get $m = ll (\Delta_q(A)) \leq
 ll(\rad^{m-1} T_q(A)) = ll (T_q(A) )-(m-1)$ and our desired inequality is proved.
 \end{pf}
 
 We shall end this subsection with some results needed later. Recall that $A^+ = C + l\rho$.
 
 \begin{lemma} Let $s \in S_l$ and
 suppose $(l-1)\rho$ belongs to the closure of the $s$-wall of $A^+$. Then $\theta_s T_q(A^+) = T_q(A^+) \oplus T_q(A^+)$.
 \end{lemma}
 
 \begin{pf} Easy weight considerations show that $T_q(A^+) = T_{-\rho}^\lambda St_q$ 
 \cite[5.2]{A00}
 and that the composite $\theta_s \circ T_{-\rho}^\lambda $ equals 
 $T_{-\rho}^\lambda \oplus T_{-\rho}^\lambda$. The lemma follows.
 \end{pf}
 
 \begin{lemma} If $s$ is as in Lemma 3.11 and $A \in \mathcal A^+$ then \\
 $\Hom_{U_q} (\theta_s L_q(A), T_q(A^+)) \simeq \begin{cases} \C^2 & \text { if } A = A^-;\\ 0 &
 \text { otherwise.} \end{cases}$
 \end{lemma}
 
 \begin{pf} By adjointness we have 
 \[
 \Hom_{U_q}(\theta_s L_q(A), T_q(A^+))\simeq \Hom_{U_q}(L_q(A),\theta_s T_q(A^+)).
 \]
  Hence the lemma follows from Lemma 3.11 via the
 fact that $\soc T_q(A^+) = L_q(A^-)$, see Proposition 3.1.
 \end{pf}
 
 We consider now a sequence of alcoves $A_0 = A^+ > A_1 > \cdots  > A_N = A^-$ as in Proposition 3.6 and let $s_i$ denote the corresponding reflections.
 \begin{prop}
 $T_q(A^+)$ is a submodule of  $\theta_{s_1} \circ \cdots \circ \theta_{s_N} L_q(A^-)$.
 \end{prop}
 \begin{pf}
 The highest weight of the tilting module $\theta_{s_1} \circ \cdots \circ \theta_{s_N} T_q(A^-)$ is $\lambda_{A^+}$ and therefore it contains $T_q(A^+)$ as a
 summand. The surjection $T_q(A^-) \rightarrow \nabla_q(A^-)$ gives a surjection $\theta_{s_1} \circ \cdots \circ \theta_{s_N} T_q(A^-) \rightarrow 
 \theta_{s_1} \circ \cdots \circ \theta_{s_N} \nabla_q(A^-)$ whose kernel $K$ does not have $L_q(A^-)$ as a composition factor. Remembering that $L_q(A^-)$ is the
 socle of $T_q(A^+)$ we conclude that the composite $T_q(A^+) \subset \theta_{s_1} \circ \cdots \circ \theta_{s_N} T_q(A^-) \rightarrow 
 \theta_{s_1} \circ \cdots \circ \theta_{s_N} \nabla_q(A^-)$ is an inclusion. If we compose this map by the surjection 
 $\theta_{s_1} \circ \cdots \circ \theta_{s_N} \nabla_q(A^-) \rightarrow \theta_{s_1} \circ \cdots \circ \theta_{s_N} (\nabla_q(A^-)/L_q(A^-))$ we get the zero map
 because $L_q(A^-)$ is not a composition factor of $\nabla_q(A^-)/L_q(A^-)$ so that by the dual version of
 Lemma 3.12 
\[
\Hom_{U_q}(T_q(A^+), \theta_{s_1} \circ \cdots \circ \theta_{s_N} (\nabla_q(A^-)/L_q(A^-))) = 0.
\]
 This means that $T_q(A^+)$ is embedded into 
 $\theta_{s_1} \circ \cdots \circ \theta_{s_N} L_q(A^-)$ as desired.
 \end{pf}
 
 \begin{remark}
 Recall that $T_q(A^+)$ is injective in $\cc_q$. Hence the statement in Proposition 3.13 says that $T_q(A^+)$ is a direct summand of 
 $\theta_{s_1} \circ \cdots \circ \theta_{s_N} L_q(A^-)$.
 \end{remark}

 \subsection{Kazhdan-Lusztig theory and Loewy length}
 We preserve the notation from above. Recall that the Lusztig conjecture on the irreducible characters, $\{\ch L_q(A)\}_{A \in \mathcal A^+}$ is known, see
 \cite {KL} or \cite {ABG} or \cite{F1}. Equivalently, the Vogan conjecture holds, i.e., we have as in
 \cite{A86}
 \begin{thm}
 If $A$ is an alcove in $X^+$ and $s \in S_l$ satisfies $As > A$ then $ll (\theta_sL_q(A)) = 3$.
 \end{thm}
 \begin{cor}
 If $\lambda \in C$ and $M \in \mathcal B_\lambda,\; s \in S_l$ then $ll(\theta_s M) \leq ll (M )+ 2.$
 \end{cor}
 
 \begin{pf} \cite{I}/\cite[D.2]{RAG}.
 \end{pf}
 
 As before $N$ denotes the number of positive roots. 
 \begin{prop}
 $ll (T_q(A^+)) = 2N + 1$ and $ll (\nabla_q(A^+)) = N + 1.$
 \end{prop}
 
 \begin{pf}
 Corollary 3.15 combined with Proposition 3.13 give immediately the inequality $ll (T_q(A^+)) \leq 2N + 1$. Then Proposition 3.10 and Corollary 3.9 give equality.
 Applying Proposition 3.10 again then implies  $ll (\nabla_q(A^+)) \leq N + 1$ and Corollary 3.9 shows that indeed we have equality also here.
 \end{pf}
 
 \begin{thm}
 For any alcove $A \in \mathcal A^{++}$ we have $ll (T_q(A)) = 2N + 1$ and $ll (\nabla_q(A)) = N+1$.
 \end{thm}
 
 \begin{pf}
 We shall proceed by induction on $A$. The induction start is provided by Proposition 3.16. So assume now the result holds for $A$ and choose an $s$-wall of $A$ with
 $A<As \subset
 (\l-1)\rho+X^+$. Weight considerations show that $T_q(As)$ is then a summand of $\theta_s T_q(A)$. Moreover, the socle of $T_q(A)$ as well as the head is $L_q(B)$ for some
 alcove $B \subset X^+$ with $sB < B$, see Proposition 3.1. Then $\theta_s T_q(A) = \theta_s(\rad T_q(A)/\soc T_q(A))$ and hence Corollary 3.15  implies that 
 $ll (T_q(As)) \leq ll (\theta_s T_q(A))) \leq ll (\theta_s(\rad T_q(A)/\soc T_q(A)) \leq 2N - 1 + 2$. Proposition 3.10 then gives $ll (\nabla_q(As)) \leq N + 1$.
 
 Suppose this last inequality was strict. Then Corollary 3.15 implies $ll (\theta_s(\nabla_q (As))) = ll (\theta_s (\nabla_q(As)/L_q(As))) \leq N + 1$. On the other 
 hand,
 $\nabla_q(A)$ is a submodule of $\theta_s \nabla_q(As)$ (see Proposition 3.4) and by induction hypothesis $ll (\nabla_q(A)) = N + 1$. Hence the head $L_q(B)$ of 
 $\nabla_q(A)$ is contained in the head of $\theta_s \nabla_q(As)$ contradicting the fact that \begin{align*}
 \Hom_{U_q}(\theta_s \nabla_q(As), L_q(B)) &\simeq
 \Hom_{U_q}(\nabla_q(As), \theta_s L_q(B)) 
 \\
 &= \Hom_{U_q}(\nabla_q(As), 0) = 0.
 \end{align*} 
 Now both equalities follow (Proposition 3.10)
 \end{pf}

  \begin{cor}
$\forall A\in\cA^{++}$
$\forall s\in S_l$,
\[
ll(\theta_s\nabla_q(A))=N+2.
\]
 \end{cor}
 
 \begin{pf}
Assume $As>A$.
Then
\begin{align*}
ll(\theta_s\nabla_q(A))
&=
ll(\theta_s\nabla_q(As))
=
ll(\theta_s(\nabla_q(As)/\soc\nabla_q(As)))
\\
&\leq
ll(\nabla_q(As)/\soc\nabla_q(As))+2
\quad\text{by Corollary 3.15}
\\
&=
N+1-1+2
\quad\text{by Theorem 3.17}
\\
&=
N+2.
\end{align*} 
On the other hand,
$\nabla_q(A)\subset \theta_s\nabla_q(As)$
by
Proposition 3.4.iii).
It follows from Theorem 3.17 again that
\[
ll(\theta_s\nabla_q(As))\geq
ll(\nabla_q(A))=N+1.
\]
Just suppose 
$ll(\theta_s\nabla_q(As))=N+1$.
Then
$\nabla_q(A)\not\subset \rad\theta_s\nabla_q(As)$.
As
$
L_q(\tilde A)
=
\hd\nabla_q(A)
\subset
\hd\theta_s\nabla_q(As)$
for some $\tilde A$
with
$\tilde As<\tilde A$,
\[
0
\not=
\Hom_{\cC_q}(\theta_s\nabla_q(As),
L_q(\tilde A))
\simeq
\Hom_{\cC_q}(\nabla_q(As),
\theta_s L_q(\tilde A))
=0,
\]
and we have a contradiction.

Assume next $As<A$.
Then 
\begin{align*}
ll(\theta_s\nabla_q(A))
&=
ll(\theta_s(\nabla_q(A)/\soc\nabla_q(A)))
\\
&\leq
ll(\nabla_q(A)/\soc\nabla_q(A))+2
\quad\text{by Corollary 3.15}
\\
&=
N+2
\quad\text{by Theorem 3.17}.
\end{align*} 
On the other hand,
\[
ll(\theta_s\nabla_q(A))
\geq
ll(\nabla_q(A))=N+1
\quad\text{
as $\theta_s\nabla_q(A)\twoheadrightarrow
\nabla_q(A)$}.
\]
Just suppose 
$ll(\theta_s\nabla_q(A))=N+1$.
Then
$\soc\theta_s\nabla_q(A)\not\subset 
\nabla_q(As)$.
It follows that
$\{\soc\theta_s\nabla_q(A)+\nabla_q(As)\}/
\nabla_q(As)\simeq
L_q(A)$,
and hence
$L_q(A)\subset \soc\theta_s\nabla_q(A)$.
But then we get a contradiction because
\begin{align*}
0
&\not=
\Hom_{\cC_q}(L_q(A),
\theta_s\nabla_q(A))
\simeq
\Hom_{\cC_q}(\theta_s L_q(A),
\nabla_q(A))
\\
&=0
\quad\text{as $As< A$}.
\end{align*}

\end{pf}

 \begin{cor}
 Let $M$ be any module in $\cc_q$. Then  $ll M \leq 2N + 1$.
 \end{cor}
 
 \begin{pf}
 Since $\cc_q$ has enough projectives (Proposition 3.1) it is enough to prove the corollary when $M$ is an indecomposable projective module. If $M \in
 \mathcal B_\lambda$ then this claim is contained in Theorem 3.17. Now for an arbitrary $\nu$ the projectives in $\mathcal B_\nu$ are summands of $T_\lambda^\nu M$ for
 appropriate $M \in \mathcal B_\lambda$. By Proposition 3.4 we have $ll (T_\lambda^\nu M) \leq ll (M)$.
 \end{pf} 
 
\subsection{Rigidity}

We shall now prove that indecomposable projective tilting modules in $\cC_q$ are rigid. The arguments are similar to those used in \cite{AK} and we 
shall therefore only give the main points.

\begin{lemma}
Let $A \in \cA^+$ and let $\nu = (l-1)\rho + l\mu$ be a special point in $X^+$. Set $A_\nu^- = A^- + l\mu$ and denote by $W_\nu$ the stabilizer of $\nu$ in $W_l$. 
If $\mu \in \rho + X^+$ then 
\begin{multline*}
 [\soc_j \nabla_q(A): L_q(A_\nu^-)] = [\rad_{N+1-j} \nabla_q(A): L_q(A_\nu^-)] \\ = \begin{cases} 1 & \text {if there exists } w \in W_\nu \text { with } l(w) = j-1 \text { and } A = w\cdot
A_\nu^-, \\
0 & \text {otherwise} .
\end{cases}
\end{multline*}
\end{lemma}

\begin{pf} 
Observe first that according to Proposition 3.4 translation onto the special point $\nu$ takes $L_q(A_\nu^-)$ to $L_q(\nu)$. The same proposition then 
shows that $[\nabla_q(A) : L_q(A_\nu^-)] = 0$ unless $\nu \in \bar A$, i.e unless $A = A_i$ for some $A_i$ occurring in a sequence of alcoves 
$\{A_i\}_{i= 1, \cdots , N}$ like in Proposition 3.7. Moreover, we get $[\nabla_q(A_i) : L_q(A^-)] = 1$ for all $i$. By Theorem 3.17 we have 
$ll \nabla_q (A_0) = N + 1$ and hence the 
arguments in the proof of Proposition 3.7 give 
$ll \I (\varphi_i \circ \varphi_{i-1} \circ \cdots \circ \varphi_1) = N+ 1 - i$ for all $i$. This means that 
$  \I (\varphi_i \circ \varphi_{i-1} \circ \cdots \circ \varphi_1) \subset \soc^{N+1-i} \nabla (A_i)$ whereas 
$  \I (\varphi_i \circ \varphi_{i-1} \circ \cdots \circ \varphi_1) \not \subset \soc^{N-i} \nabla (A_i)$. Since $\hd \nabla_q(A_0) = L_q(A_\nu^-)$ 
this proves the
statement about the occurence of $L_q(A_\nu^-)$ in the socle layers of $\nabla_q(A_i)$. To get the statement about the radical series of $\nabla_q(A_i)$ we 
argue similarly using this time the composite 
$\varphi_N \circ \varphi_{N-1} \circ \cdots \circ \varphi_{i+1} : \nabla_q(A_i) \rightarrow \nabla_q(A_\nu^-)$. In fact, 
if
$L_q(A_\nu^-)\subset
\rad_{j}\nabla_q(A_{i})$,
then $\varphi_N\circ\dots\circ\varphi_{i+1}(\rad^j\nabla_q(A_{i})) \neq 0$. Moreover, 
the assumption $\mu \in \rho + X^+$ ensures via Theorem 3.17 that 
$ll \nabla(A_i) = N+1 $ for all $i$. Since each $\varphi_i$ kills the socle we get 
\[
1 \leq ll(\varphi_N\circ\dots\circ\varphi_{i+1}(\rad^j\nabla_q(A_{i})))
\leq
ll(\rad^j\nabla_q(A_{i}))-(N-i) = i+1-j. \]
Hence $j \leq i$ for all such $j$ and the statement follows. 

\end{pf}

 \begin{defn}
 Let $M \in \cC_q$. We say that a composition factor $L$ of $M$ is rigidly placed in $M$ if we have
 $ [\rad^i M : L] = [\soc^{ll M - i} : L] $ for all $i$. 
 \end{defn}
 
 Clearly a module $M$ is rigid iff all its composition factors are rigidly placed.

 \begin{prop}
Let $s \in S_l$ and assume
$A, As\in\cA^{++}$
and that 
$L$ is a simple module in $\cC_q$.
If $L$ is rigidly placed in
both $\nabla_q(A)$
and $\nabla_q(As)$,
then
$L$ is also rigidly placed in $\theta_s\nabla_q(A)$.

\end{prop}

\begin{pf}
We may assume that $As < A$. 
Recall from
Proposition 3.4 that we then have a short exact sequence
\[
0\to
\nabla_q(As)
\to
\theta_s\nabla_q(A)
\overset{\pi}{\to}
\nabla_q(A)\to0.
\]
As observed in the proofs of Theorem 3.17 and Corollary 3.18
we have 
$\hd \theta_s\nabla_q(A) = \hd \nabla_q (A)$ and $\soc \theta_s\nabla_q(A) = \soc \nabla_q (As)$. This implies

\begin{align*}
\rad^{i-1}\nabla_q(As)
&\subset
\rad^{i-1}(\rad\theta_s\nabla_q(A))
=
\rad^i\theta_s\nabla_q(A),
\\
\text{and } 
\pi(\soc^{N+2-i}\theta_s\nabla_q(A))
&\subset \soc^{N+1-i}\nabla_q(A).
\end{align*}
On the other hand,
\begin{align*}
\nabla_q(As)\cap\soc^{N+2-i}\theta_s\nabla_q(A)
&=
\soc^{N+2-i}\nabla_q(As),
\\
\pi(\rad^i\theta_s\nabla_q(A))
&=
\rad^i\nabla_q(A).
\end{align*}
It follows that
\begin{align*}
[\soc^{N+2-i}
&\theta_s\nabla_q(A) 
:
L]
\leq
[\soc^{N+2-i}\nabla_q(As) :
L]
+
[\soc^{N+1-i}\nabla_q(A) :
L]
\\
&=
[\rad^{i-1}\nabla_q(As) :
L]
+
[\rad^i\nabla_q(A) :
L]
\quad\text{by the hypothesis}
\\
&\leq
[\rad^{i}\theta_s\nabla_q(A) :
L].
\end{align*}
By Corollary 3.18 we have $ll \theta_s \nabla_q((A) = N+2$ and hence this inequality 
proves that $L$  is rigidly placed in
$\theta_s\nabla_q(A)$.
\end{pf}

\begin{prop}
Let $s \in S_l$ and 
$A\in\cA^{++}$. Suppose
$B, Bs \in\cA^+$ with $B > Bs$.
If $L_q(B)$ is rigidly placed in
$\theta_s\nabla_q(A)$,
then
$L_q(Bs)$ is rigidly placed 
in
$\nabla_q(A)$.

\end{prop}

\begin{pf}
We assume first that $A < As$ so that $\nabla_q(A) \subset \theta_s \nabla_q(A)$. We set 
$L=L_q(B)$.
As $Bs<B$,
$\theta_sL=0$. This means in particular that $L$ cannot occur in the head or socle
of $\theta_s V$ for any $V \in \cC_q$.

Write $M=\nabla_q(A)$.
Then we have for all $i\in\N$,
\begin{align*}
ll(\theta_s\soc^{N+1-i}M)
&\leq
ll(\soc^{N+1-i}M)+2
\quad\text{by Corollary 3.15}
\\
&\leq
N+3-i,
\end{align*}
and hence
$\theta_s\soc^{N+1-i}M\subset
\soc^{N+3-i}\theta_sM$.
This implies
\begin{align*}
\rad\theta_s\soc^{N+1-i}M
&\subset
\soc^{ll(\theta_s\soc^{N+1-i}M)-1}
\theta_s\soc^{N+1-i}M
\\
&\subset
\soc^{N+2-i}
\theta_s\soc^{N+1-i}M
\\
&\subset
\soc^{N+2-i}
(\soc^{N+3-i}\theta_sM)
=
\soc^{N+2-i}
\theta_sM.
\end{align*}
Since (as observed above) $L$ does not occur in the head of $\theta_s\soc^{N+1-i}M$ we get from this
\[
[\theta_s\soc^{N+1-i}M : L]
\leq
[\soc^{N+2-i}\theta_sM : L].
\tag{1}
\]
Let now $M_i$ denote the submodule of $\theta_s M$ containing $\theta_s \rad^i M$
with $M_i/\theta_s\rad^iM = \soc(\theta_sM/\theta_s\rad^iM)$.
Then
\begin{align*}
ll(\theta_sM/\theta_s\rad^iM)
&=
ll(\theta_s\rad M/\theta_s\rad^iM)
\quad\text{as $\theta_sM=\theta_s\rad M$}
\\
&=
ll(\theta_s(\rad M/\rad^iM))
\\
&\leq
ll(\rad M/\rad^iM)+2
\quad\text{by Corollary 3.15}
\\
&=
i-1+2=i+1.
\end{align*}
It follows that
\[
\rad^i(\theta_sM/\theta_s\rad^iM)
\subset 
\soc(\theta_sM/\theta_s\rad^iM)
=
M_i/\theta_s\rad^iM,
\]
and hence
$\rad^i\theta_sM \subset M_i$.
As $L$ does not occur in 
$M_i/\theta_s\rad^iM =\soc(\theta_s(M/\rad^iM)$
we find 
\[
[\rad^i\theta_sM : L]\leq
[M_i : L]=
[\theta_s\rad^iM : L].
\tag{2}
\]
Then
\begin{align*}
[\theta_s\soc^{N+1-i}M : L]
&\leq
[\soc^{N+2-i}\theta_sM : L]
\quad\text{by (1)}
\\
&=
[\rad^i\theta_sM : L]
\quad\text{as
$L$ is rigidly placed in
$\theta_sM$}
\\
&\quad\text{and as
$ll(\theta_sM)=N+2$
by
Corollary 3.18}
\\
&\leq
[\theta_s\rad^iM : L]
\quad\text{by (2)}
\\
&\leq
[\theta_s\soc^{N+1-i}M : L],
\end{align*}
i.e.,
\begin{multline*}
[\theta_s\soc^{N+1-i}M : L]
=
[\soc^{N+2-i}\theta_sM : L]
=
\\
[\rad^{i}\theta_sM : L]
=
[\theta_s\rad^iM : L].
\end{multline*}
In particular,
$
[\theta_s(\soc^{N+1-i}M/\rad^iM) : L]=0$.
But this means
\[
\sum_{D\in\cA^+}
[\soc^{N+1-i}M/\rad^iM : L_q(D)]
[\theta_sL_q(D) : L] = 0,
\]
and we conclude that if for some $D$ we have $[\theta_sL_q(D) : L]\ne0$, then 
$[\soc^{N+1-i}M/\rad^iM : L_q(D)]=0$.
However, 
$[\theta_sL_q(Bs) : L]=1$, so we get in particular,
$
[\soc^{N+1-i}M/\rad^iM : L_q(Bs)]=0$, i.e.
$L_q(Bs)$
is rigidly placed in $M=\nabla_q(A)$.
In the case where $A > As$ so that $\nabla_q(A)$ is a quotient of $\theta_s \nabla_q(A)$ 
a completely analogous argument works (cf. \cite{AK}). 
\end{pf}

The above results now give us the following

\begin{thm}
If
$A$
is sufficiently deep inside $\cA^{+}$ then $\nabla_q(A)$
is rigid.
\end{thm}

\begin{pf}
Let $L = L_q(B)$ for some $B \in \cA^+$. Choose a special point $\nu \in X^+$ such that $B \subset \nu - X_l$. 
The top alcove in this box is $A_\nu^-$ and in general we can find a sequence of alcoves $B_0 = B < Bs_1 < \cdots Bs_1s_2\cdots s_m = A_\nu^-$ 
with $s_i \in S_l$
inside $\nu - X_l$. We want to prove that $L$ is rigidly placed in any $\nabla_q(A)$ with
$A$ sufficiently deep inside $\cA^+$. Suppose $A$ is such an alcove for which $L$ is not rigidly placed in $\nabla_q(A)$. Since $A$ is 
sufficiently deep inside $\cA^{++}$ we have certainly $As_1 \in \cA^{++}$ and more generally $As_{j_1} \cdots s_{j_r} \in \cA^{++}$ for any subsequence $s_{j_1},
\cdots, s_{j_r}$ of $s_1, s_2, \dots, s_m$. Hence
Propositions 3.22 and 3.23 imply that $L_q(Bs_1)$ is not rigidly placed either in $\nabla_q(A)$ or in $\nabla_q(As_1)$ and repeating
this argument $m$ times that there exists a subsequence $s_{j_1}, \cdots, s_{j_r}$ of $s_1, s_2, \cdots, s_m$ such that $L_q(A_\nu^-)$ is not
rigidly placed in $\nabla_q(As_{j_1} \cdots s_{j_r})$. But then $ As_{j_1} \cdots s_{j_r} = w\cdot A_\nu^-$ (otherwise $L_q(A_\nu^-)$ is not a composition
factor of $\nabla_q(As_{j_1} \cdots s_{j_r})$) for some $w \in W_\nu$ and we have a 
contradiction to Lemma 3.20.

\end{pf}
 
\begin{rem}
The wording `sufficiently deep' in Theorem 3.24  means that $A \subset n\rho + X^+$ for some large integer $n$. Actually we 
believe $n = l$ should be enough (i.e. $A \in \cA^{++}$) but (cf. the proof below) our present arguments 
require $n = hl$. On the other hand, our examples in Section 5 
demonstrate that we cannot take $n$ smaller than $l$.

To see that $n=hl$ suffices we need the following fact. We use the notation and assumptions from the proof of Theorem 3.24. In particular,
$B$ is an alcove inside $\nu - X_l$ and $s_1,s_2,\cdots, s_m$ is a minimal sequence such that 
$Bs_1 s_2 \cdots s_m = A^-_\nu$.
In particular,
$B$ is an alcove inside $\nu - X_l$ and $s_1,s_2,\cdots, s_m$ is a minimal sequence such that 
$Bs_1 s_2 \cdots s_m = A^-_\nu$.
Let
$A \in hl\rho + \cA^+$. 
Then any of the subsequence $s_{j_1}, \cdots, s_{j_r}$ of $s_1, s_2, \cdots, s_m$ we have $As_{j_1} \cdots s_{j_r} \in \cA^{++}$.

\begin{pf} 
Recall that $A^-_\nu=w_0\cdot C+\nu + \rho$ and set  
$A^+_{\nu 
-l\rho}=C+\nu+\rho-l\rho$. We can write
$A=w\cdot
A^+_{\nu-l\rho}+lw\gamma$
for some
$w\in W,\gamma\in\Z R$.
Then
$A=w\cdot
C+w(\nu-l\rho+\rho+l\gamma)$,
so $w(\nu-l\rho+\rho+l\gamma) \in
l(h\rho+X^+)$.
Also $As_{j_1}\dots s_{j_r}
=
w\cdot
A^+_{\nu-l\rho}s_{j_1}\dots s_{j_r}+
wl\gamma$.
We claim
$A^+_{\nu-l\rho}s_{j_1}\dots s_{j_r}\in
\nu + \rho -l\rho+W\cdot\{
A'\in\cA^+\mid 
A'\leq
A^-_{l\rho-\rho}\}$.
Indeed, if for some $1 \leq k \leq r$ we have
$A^+_{\nu-l\rho}s_{j_1}\dots s_{j_k-1}\in
\nu + \rho -l\rho+\cA^+$
with
$A^+_{\nu  -l\rho}s_{j_1}\dots s_{j_k-1}s_{j_k}
\not\in
\nu+\rho-l\rho+\cA^+$, then
$A^+_{\nu-l\rho}s_{j_1}\dots s_{j_k-1}s_{j_k}
=
\nu +\rho -l\rho+s_k'\cdot
(A^+_{\nu-l\rho}s_{j_1}\dots s_{j_k-1}-\nu-l\rho-\rho)$
for some reflexion $s_{j_k}'\in W$ in a hyperplane
orthogonal to a simple root.
Hence we may assume
$A^+_{\nu-l\rho}s_{j_1}\dots s_{j_r}\in
\nu+\rho-l\rho+\cA^+$.
Then
by \cite[1.6]{V},
$A^+_{\nu-l\rho}s_{j_1}\dots s_{j_r}\leq
A^-_{\nu}$
as
$\rd(A^+_{\nu-l\rho}, B)+
\rd(B,A^-_\nu)=
\rd(A^+_{\nu-l\rho},A^-_\nu)$.
Thus we can write
\[
As_{j_1}\dots s_{j_r}
=
w\cdot
(w'\cdot
A'+{\nu-l\rho}+\rho)+
wl\gamma
=
ww'\cdot
A'+w(\nu-l\rho+\rho+l\gamma)
\]
for some
$w'\in W$ and $A'\in\cA^+$ with
$A'\leq A^-_{l\rho-\rho}$.
It follows for any simple root $\alpha$
that
\begin{align*}
\langle
ww'\cdot
\lambda_{A'}
&+w(\nu-l\rho+\rho+l\gamma)+\rho,
\alpha^\vee\rangle
\\
&=
\langle
\lambda_{A'}+\rho,
(ww')^{-1},\alpha^\vee\rangle
+\langle
w(\nu-l\rho+\rho+l\gamma),
\alpha^\vee\rangle
\\
&\geq
-\langle
l\rho,
\alpha_0^\vee\rangle
+
lh
=
l,
\end{align*}
as desired.

\end{pf}

\end{rem}

 Finally we can formulate our main result on rigidity of certain indecomposable tilting modules.
 Recall that there are examples of non-rigid tilting modules, see Section 5. So at least for the regular block our result  
 here seems to be the best possible.
 
 \begin{thm}
 Let $A \in \cA^{++}$. Then $T_q(A)$ is rigid.
 \end{thm}
 
 \begin{pf}
 When $A$ is sufficiently deep inside $\cA^+$ the result follows from Theorem 3.24 exactly as in \cite {AK}. So here we shall
 only show how to deduce the general case from this.
 
 Recall that for any $\mu \in X^+$ we have \cite{Lu}
 \[
L_q(\mu) \simeq L_q(\mu^0) \otimes L_\C(\mu^1)^{[1]} \tag{1}
\]
where $L_\C(\mu^1)$ denotes the irreducible module with highest weight $\mu^1$ for the simple complex Lie algebra corresponding to $R$ and an upper ${[1]}$ 
  on a module means that we consider it as a $U_q$-module via the quantum Frobenius homomorphism.
  
  Likewise we have a factorization of projective tilting modules. Namely, 
 if $\eta \in (l-1)\rho + X_l$ and $\eta' \in X^+$ then we have \cite{A92}
 \[
  T_q(\eta + l\eta') \simeq T_q(\eta) \otimes L_\C(\eta')^{[1]}.
  \tag{2}
  \]
  
  Let now $A \in \cA^{++}$. Then we may write $A = A' + l\mu$ for some alcove $A'$ contained in $l\rho +X_l$ and some $\mu \in X^+$.
  Consider now $m \gg 0$. Then $L_\C(\mu) \otimes L_\C(m\rho) = \bigoplus _\eta L_\C(\eta +m\rho)$ where the sum runs over the
  multiset of weights of $L_\C(\mu)$. When $m$ is big enough the alcoves of the form $A' + l\eta + lm\rho$ all lie sufficiently deep inside $\cA^+$
  for all these $\eta$'s. Hence all $T_q(A'+ l\eta + lm\rho)$ are rigid. But by (2) we see that $T_q(A) \otimes L_\C(m\rho)^{[1]}$ equals the direct sum of these and
  since they all have the same Loewy length (namely $2N+1$ according to Theorem 3.17) this sum is also rigid. However, (1) implies that tensoring by
  any  $L_\C(\mu)^{[1]}$ takes a semisimple module in $\cC_q$ into a semisimple module. Hence we have 
 \begin{multline*} (\rad^i(T_q(A) \otimes L_\C(m\rho)^{[1]}) \subset  (\rad^i T_q(A)) \otimes L_\C(m\rho)^{[1]} \subset \\ (\soc^{2N+1-i} T_q(A)) \otimes L_\C(m\rho)^{[1]}) 
 \subset  \soc^{2N+1-i} (T_q(A) \otimes L_\C(m\rho)^{[1]})
 \end{multline*}
 By the above the two outer terms are equal so that we must have equality everywhere. We conclude that $\rad^i T_q(A) =\soc^{2N+1-i} T_q(A)$, i.e., $T_q(A)$ is rigid.
 
 \end{pf}

\begin{remark}
To obtain Theorem 3.26 only,
we could make use of the fact from
\cite[3.4]{APW92}
that
$
\soc_{\cC_q}M=\soc_{\bar\cC_q}M$
for each
$M\in\cC_q$,
where
$\bar\cC_q$ 
is the category of finite
dimensional modules of type 1 over the small
quantum algebra, see 5.2.

\end{remark}

  \section{Loewy structures in the modular case}

In this section we shall consider the modular category $\cC$ introduced in Section 2.2. 
We shall prove that certain tilting modules are also rigid in this case. Although we could proceed 
just as in the quantum case treated in Section 3 we shall make some shortcuts by 
taking advantage of our results in \cite {AK}. 

In addition to the notation already introduced in Section 2 we shall need some more notation.
\subsection{More notation}
We let $F: G \to G$ denote the Frobenius morphism coming from the $p$-th power map on $k$. 
The kernel of $F$ is denoted $G_1$. This is an infinitesimal subgroup scheme of $G$. We shall also consider the 
subgroup scheme $G_1T$ and its representations. 

If $M \in \cC$ we denote by $M^{(1)}$ the Frobenius twist of $M$. This is the same vector space but the $G$-action is given by $g \cdot v=
F(g) v, g \in G, v \in V$. Clearly $G_1$ acts trivially on $M^{(1)}$. If on the other hand $G_1$ acts trivially on some 
$V \in \cC$ then $V = M^{(1)}$ for some $M \in \cC$. In this case we write also $M = V^{(-1)}$. 

We shall carry 
over much of the notation from Section 3 by just replacing $l$ by $p$. In particular, $W_p$ is now the affine Weyl group in question. The
lowest alcove in $X^+$ is still denoted $C$ but now given by $C =  
\{ \lambda \in X^+ \mid \langle \lambda + \rho, \alpha^\vee \rangle < p$ for all $\alpha \in R^+ \}$. The fundamental box is denoted $X_p$.

\subsection{Rigidity for injective $G_1T$-modules}

The Lusztig conjecture for irreducible $G$-modules  may just as well be formulated as a conjecture about irreducible $G_1T$-modules. It is then 
natural to expect the conjecture to be true for $p \geq h$ but it is only known  to hold in general for $p \gg 0$, see \cite{AJS} and \cite {F1}. 
On the other hand, the example in \cite{A94} shows
that it cannot be expected to hold for $p < h$. In this section we shall assume that the conjecture does hold for our group. As in Section 3 we 
will use the analogue in form of the Vogan conjecture which we can formulate by saying that 
\[
ll \theta_s L_1(\lambda) \leq 3 \text { for all $p$-regular }
\lambda \in X \text { and all $s \in S_p$.}
\]
Here $L_1(\lambda)$ denotes the irreducible $G_1T$-module with highest weight $\lambda$. If $\lambda \in X_p$ then $L_1(\lambda)$ is just the restriction to
$G_1T$ of the irreducible $G$-module $L(\lambda)$. In general, we write $\lambda = \lambda^0 + p \lambda^1$ (in analogy with Section 3.1) with $\lambda^0 \in X_p$ and $\lambda^1
\in X$ and have $L_1(\lambda) = L(\lambda^0) \otimes p\lambda^1$.

As in Section 3 we get then
\begin{cor}
Let $s \in S_p$ and suppose $M$ is a $G_1T$-module belonging to a regular block. Then $ll \theta_s M \leq ll M + 2$.
\end{cor}
 
Denote now for any $\lambda \in X$ by $Q_1(\lambda)$ the injective envelope of $L_1(\lambda)$ (in the category of finite dimensional rational $G_1T$-modules). 
Then we have
\begin{thm} \cite[Theorem 7.2]{AK}
Let $\lambda \in X$ be $p$-regular. Then the injective $G_1T$-module $ Q_1(\lambda)$ is rigid with Loewy length $2N+1$.
\end{thm}

\subsection{Rigidity of some tilting modules for $G$}

The injective $G_1T$-modules are expected to have $G$-structures. In fact, Donkin has conjectured
\begin{conj} (S. Donkin) \cite{Do93}
Let $\lambda \in X_p$ and write $\tilde \lambda = 2(p-1)\rho + w_0\lambda$. Then $Q_1(\lambda) \simeq T(\tilde \lambda)_{\mid G_1T}$.
\end{conj}

This conjecture is known to be true for $p \geq 2(h-1)$ \cite{RAG}. In particular, it holds  under the
stronger assumption in the following corollary which we here have formulated as a statement about tilting modules using the above relation

\begin{cor} \cite[Proposition 8.4]{AK}
Suppose $p \geq 3h-3$ and let $\lambda \in X_1$ be $p$-regular. Then the $G$-module $T(\tilde \lambda)$ is rigid with Loewy length
$2N+1$. In fact, we have $\soc^j_G T(\tilde \lambda) = \soc^j_{G_1} Q_1(\lambda)$ for all $j$.
\end{cor}
We have added a subscript to the socle series to indicate which category we consider.

\begin{rem}
The assumption $p \geq 3h-3$ ensures that for all $\lambda \in X_1$ the $G$-composition factors $L(\mu)$ of $T(\tilde \lambda)$ always have $\mu^1 \in \bar C$. By 
the linkage principle this implies that the isotypic components of the $G_1T$-Loewy layers of $T(\tilde \lambda)$ are all semisimple
as $G$-modules. 
\end{rem}  

In order to 
prove rigidity for a wider class of tilting modules we need

\begin{lemma}
Let $M$ and $L$ be $G$-modules which satisfy
\begin{enumerate}
\item $\soc_G M = \soc_{G_1} M$,
and
\item $(\soc_G M) \otimes L^{(1)}$ is semisimple.
\end{enumerate}
Then $\soc_G (M \otimes L^{(1)}) = (\soc_G M) \otimes L^{(1)}$.
\end{lemma}

\begin{pf}
Since both sides (for the right hand side we use assumption (2)) of the desired equality are semisimple $G$-modules it is enough to prove
$$ \Hom_G(L(\lambda), \soc_G (M \otimes L^{(1)})) = \Hom_G(L(\lambda), (\soc_G M) \otimes L^{(1)})$$
for all $\lambda \in X^+$. 
Here the left hand side equals
\begin{multline*}
 \Hom_G(L(\lambda^0) \otimes L(\lambda^1)^{(1)}, M \otimes L^{(1)}) = 
 \\
 \Hom_G(L(\lambda^1), \Hom_{G_1}(L(\lambda^0), M)^{(-1)} \otimes L)
 \end{multline*}
whereas the right hand side equals
\begin{multline*}
\Hom_G(L(\lambda^0) \otimes L(\lambda^1)^{(1)}, (\soc_G M) \otimes L^{(1)}) = 
\\
\Hom_G(L(\lambda^1), \Hom_{G_1}(L(\lambda^0), \soc_G M)^{(-1)} \otimes L).
\end{multline*}
By assumption (1) $\soc_{G_1} M = \soc_G M$ and hence 
\[
\Hom_{G_1}(L(\lambda^0), \soc_G M) = \Hom_{G_1}(L(\lambda^0), M).
\]
The conclusion follows.
\end{pf}

This now allows us to prove
\begin{thm}
Assume $p \geq 3h-3$ and suppose $\lambda \in X$ is a $p$-regular weight which satisfies $p \leq \langle \lambda + \rho, \alpha^\vee \rangle \leq p(p-h+2)$ for 
all $\alpha \in R^+$. Then $T(\lambda) $ is rigid of Loewy length $2N+1$.
\end{thm}

\begin{pf} Using the lower bounds in our assumptions on $\lambda$ we see that we may write $\lambda = \tilde \mu + p \nu$ with $\mu \in  X_p$ and $\nu \in X^+$.
Now by \cite[2.1]{Do93} we have 
\[ T(\lambda)
\simeq T(\mu) \otimes T(\nu)^{(1)}.
\]

The upper bound on $\lambda$ ensures that any dominant weight $\eta$ of $T(\lambda)$ has $\eta^1 \in \bar C$. This means in particular that 
$\nu \in  \bar C$ so that by the linkage principle $T(\nu) = L(\nu)$. It also means that for any composition factor $L(\eta)$ of $T(\tilde \mu)$ we have $\eta^1 +
\nu \in \bar C$ and hence (again by the linkage principle) $L(\eta) \otimes L(\nu)^{(1)} = L(\eta^0) \otimes (L(\eta^1) \otimes L(\nu))^{(1)}$ is semisimple. Therefore
for any $j$ the pair $(T(\tilde \mu)/\soc^j_G T(\tilde \mu), T(\nu))$ satisfies condition (2) in Lemma 4.6. Moreover, by the corollary above
the $G$-socle series of $T(\tilde \mu)$ coincides with the $G_1$ socle series. Hence condition (1) in Lemma 4.6 holds for each 
$ T(\tilde \mu)/\soc^j_G T(\tilde \mu)$. Using this we see that the lemma combined with Corollary 4.4 and Theorem 4.2 give
\begin{align*}
\soc_G^{j+1} T(\lambda) 
&= (\soc^{j+1}_{G_1} Q_1(\mu)) \otimes T(\nu)^{(1)} 
\\
&= (\rad_{G_1}^{2N-j} Q_1(\mu)) \otimes T(\nu)^{(1)} = \rad_G^{2N-j} T(\lambda).
\end{align*}
Here the last equality uses the `dual version' of Lemma 2.4.
\end{pf}

\begin{rem}
In Section 6 we give an example of a non-rigid tilting module for $G = SL_3(k)$. This example shows that the upper bound in Theorem 4.7 is just about `best possible'.
\end{rem}

\section{Type $B_2$}

In this section we calculate the Loewy structure of all regular Weyl modules in $\cC_q$ in the case where $R$ is of type $B_2$. Moreover, we compute in this
situation also the Loewy structure of several indecomposable tilting modules 
with small highest weights. It turns out that there is exactly one alcove in
$\cA^+$ for which the corresponding Weyl module is not rigid. We prove that
the indecomposable tilting
module associated with this alcove is likewise non-rigid. 

We start out by two subsections containing some results which apply to the general case. For Weyl modules of small Loewy lengths the strategy here is very effective
for proving rigidity. In particular, we shall deduce that all regular Weyl modules are rigid for type $A_2$ and with the single exception mentioned above we
get the same result for type $B_2$.

\subsection{Parity filtrations} 
\begin{defn}
Let $M \in \cC_q$. We say that a filtration $0 = M^r \subset M^{r-1} \subset \cdots \subset M^0 = M$ with layers $M_i = M^i/M^{i+1}$ is a parity filtration 
if $d(A,B) + i -j$ is even whenever $L_q(A)$ , respectively $L_q(B)$, is a composition factor of $M_i$, respectively of $M_j$.
\end{defn}

\begin{rem}
Let $M^\bullet$ be a parity filtration of $M \in \cC_q$ and let $r > i \geq 0$. Then we have in particular that $d(A,B)$ is even for all composition 
factors of $M_i$. According to the $q$-analogue of Corollary 2.10 in \cite{A86} this means that $M_i$ is semisimple. So by the general properties of radical and socle 
series we have $\rad^i M \subset M^i \subset \soc^{r-i} M$. If $r = ll M$ then $M^\bullet$ is a Loewy series for $M$ and $M$ is rigid if this Loewy series
coincides with the socle and the radical filtrations.
\end{rem}

\begin{prop}
Suppose $M \in \cC_q$ belongs to a regular block. Assume that 
\begin{enumerate}
\item
$M$ has a parity filtration of length $r$
\item
$M$ has simple head and socle.
\end{enumerate} 
If $r \leq 4$ then $M$ is rigid.
\end{prop}

\begin{pf}
Denote by $L_q(A)$ the socle of $M$ and by $L_q(B)$ the head of $M$. 
Then by assumption (2) we have $\soc M = L_q(A) = M^{r-1} = \rad^{ll M - 1} M$ and 
$\rad M = M^1 = \soc^{ll M-1} M$. The fact that two simple modules $L_q(A_1)$ and $L_q(A_2)$ do not extend 
when $d(A_1, A_2)$ is even implies that the layer $\soc_2 M$ consists of composition factors $L_q(D)$ in $M$ having $d(A,D)$ odd. Similarly, all 
composition factors $L_q(D)$ of 
$\rad_1 M$ have $d(D, B)$ odd.

The proposition is trivial if $r \leq 2$. If $r = 3$, we may assume the filtration is strict, meaning
$0=M^3\subsetneq M^2\subsetneq M^1\subsetneq M^0 = M$.
Then $M_1$ consists of all composition factors $L_q(D)$ of $M$ for which $d(A,D)$ is odd. By the
above this means that $\soc_2 M = M_1$. So we have $\soc^2 M = M^1 = \rad M$, i.e. $M$ is rigid.

If $r=4$  the parity condition ensures that $M_2$, respectively $M_1$ consists of composition factors $L_q(D)$ with $d(A,D)$ odd and $d(D,B)$ even, 
respectively $d(D,B)$ odd and $d(A,D)$ even. The 
same argument as above gives then $\soc^2 M = M^2$ and $\rad^2 M = M^2$. Rigidity of $M$ follows.
\end{pf}

\subsection{Inducing from small quantum groups}
In \cite{AK} we proved rigidity of regular baby Verma modules for $G$. The techniques there carry over to the quantum case. To be precise,
let $\mathfrak{u}_q$
(resp. $U_q^0$) be the subalgebra of
$U_q$
generated by
$E_i,F_i,K_i$ 
(resp.
$K_i$
and 
$\lbr
K_i
\\
t
\rbr$),
and let
$U_q^{\leq0}$
be the subalgebra generated by
$U_q^0$
and $F_i^{(t)}$,
$i\in[1,n]$,
$t\in\N$.
Let
$\hat\cC_q$
(resp. $\tilde\cC_q$,
$\overline{\cC}_q$)
be the category of finite dimensional
$\mathfrak{u}_qU_q^{0}$-
(resp. $\mathfrak{u}_qU_q^{\leq0}$-,
$\mathfrak{u}_q$-)
modules of type $\bf1$ as in
\cite[0.4]{APW92},
and let
$\tilde\nabla_q$
(resp. $\mathrm{ind}_{\tilde\cC_q}^{\cC_q}$)
be the induction functor
from
$\cC_q^{\leq0}$ to $\tilde\cC_q$
(resp. $\tilde\cC_q$ to $\cC_q$). Then we have an induced modules (or quantized baby Verma modules) $\tilde \nabla_q(A)$ corresponding to each alcove $A \in \cA$ and 
the analogue of
Theorem 5.6 in \cite{AK} says that all these modules are rigid as $\mathfrak {u}_q U^0_q$-modules.  

Set now $\tilde \nabla_q(A)^i = \rad^i_{\hat \cC_q} \tilde \nabla_q(A)$ for $i \geq 0$. These submodules belong to $\tilde \cC_q$ and in the usual notation for
layers we get the following exact sequences
\[
 0 \to \ind_{\tilde \cC_q}^{\cC_q} (\tilde \nabla_q(A)^{i+1}) \to \ind_{\tilde \cC_q}^{\cC_q} (\tilde \nabla_q(A)^i) 
\to \ind_{\tilde \cC_q}^{\cC_q} (\tilde \nabla_q(A)_i) . \tag{1}
\]
Noting that $\ind_{\tilde \cC_q}^{\cC_q}$ takes the simple module $\tilde L_q(A)$ in $\tilde \cC_q$ to $L_q(A)$ if $A \in \cA^+$ and to $0$ otherwise, we get a filtration
$\nabla_q(A)^\bullet = \ind_{\tilde \cC_q}^{\cC_q} (\tilde \nabla_q(A)^\bullet)$ of $\nabla_q(A)$ with semisimple layers.

\begin{prop}
Let $A \in \cA^+$. Then $\nabla_q(A)$ has a parity filtration. Moreover, $ll \nabla_q(A) = N+1$ if $A \in \cA^{++}$ and $ll \nabla_q(A) < N+1$ if $A \in \cA^+ \setminus \cA^{++}$.
In particular, if
$A\in\cA^{++}$, the length of the parity filtration on $\Delta_q(A)$ equals its Loewy length.
\end{prop}

\begin{pf}  The $q$-analogues of the results from \cite{AK} give that the 
$\hat\cC_q$-Loewy series of $\tilde \nabla_q(A)$ is a parity filtration for all $A$ . Then (1) above shows
that so is the filtration $\nabla_q(A)^\bullet$ because the composition factors $L_q(B)$ of $\nabla_q(A)_i$ are among those for which $\tilde L_q(B)$ is a
composition factor of $\tilde \nabla_q(A)_i$.

The Loewy length equality is already proved in Theorem 3.17. The inequality follows from the above by observing that the top layer $\tilde \nabla_q(A)_0 =
\hd \tilde \nabla_q(A)$ is
$\tilde L_q(\tilde A)$ and when $A \in \cA^+ \setminus \cA^{++}$ the alcove $\tilde A$ lies outside $\cA^+$ so that $\nabla_q(A) = \nabla_q(A)^0=\nabla_q(A)^1$.
\end{pf}

In type $A_2$ and $B_2$
the $\tilde\cC_q$-structure on each
$\tilde\nabla_q(A)_i$
from \cite{KY07}/\cite{KY}
shows that
the map
$\ind_{\tilde \cC_q}^{\cC_q} (\tilde \nabla_q(A)^i) 
\to \ind_{\tilde \cC_q}^{\cC_q} (\tilde \nabla_q(A)_i)$
of (1)
is surjective and that
$\ind_{\tilde \cC_q}^{\cC_q} (\tilde \nabla_q(A)_i)$
is explicitly obtained by the Borel-Weil-Bott
cancellation.

\begin{cor}
Suppose $R$ is of type $A_2$. Then $\nabla_q(A)$ is rigid for all $A \in \cA^+$.
\end{cor}

\begin{pf} Let $M = \nabla_q(A)$ with $A \in \cA^+$. Then Proposition 5.4 says that assumption (1) in Proposition 5.3 holds for $M$. Moreover, it is well known that
Weyl modules for Type $A_2$ have simple socles (cf. Appendix where we deduce this from a general proposition) so that also assumption (2) holds. Hence 
Proposition 5.3 gives the result.  

\end{pf}

\begin{cor} Suppose $R$ is of type $B_2$. Let $A \in \cA^+$ be an alcove for which $\nabla_q(A)$ has simple head. Then $\nabla_q(A)$ is rigid.
\end{cor}

\begin{pf} Suppose first that $A \in \cA^+ \setminus \cA^{++}$. Then Proposition 5.4 implies that $ll \nabla_q(A) \leq 4$ and just as in Corollary 5.5 we can 
apply Proposition 5.3.

Let next $A \in \cA^{++}$. Then Theorem 3.17 gives $ll \nabla_q(A) = 5$ and we are out of the range where Proposition 5.3 applies. However, we can use a variation of
the arguments: We still have a parity filtration (Proposition 5.4) with $\rad^4 \nabla_q(A) = \nabla_q(A)^4 = \soc \nabla_q(A)$ and $\rad \nabla_q(A) = \nabla_q(A)^1
= \soc^4 \nabla_q(A)$. The parity arguments from the proof of Proposition 5.3 show that $\soc_2 \nabla_q(A)$ has no composition factors in common with $\nabla_q(A)_2$
and similarly $\rad_2 \nabla_q(A)$ has no composition factors in common with $\nabla_q(A)_2$. We want $\soc^2 \nabla_q(A) = \nabla_q(A)^3 = \rad^3 \nabla_q(A)$
The first equality
follows if we check $\Ext_{\cC_q}^1(L_q(B), L_q(A)) = 0$ for all
composition factors $L_q(B)$ of $\nabla_q(A)_1$
while the second follows if $\Ext_{\cC_q}^1(L_q(B), L_q(\tilde A)) = 0$ for all
composition factors $L_q(B)$ of $\nabla_q(A)_3$. 
Now the vanishing of these $\Ext^1$'s is a consequence of the Kazhdan-Lusztig theory. In fact, 
$\Ext_{\cC_q}^1(L_q(B), L_q(A)) \subset
\Ext_{\cC_q}^1(L_q(A), \nabla_q(B))$
for
$B<A$,
and the dimension of the latter is given as the ``leading" coefficient $\mu(B, A)$
of the relevant Kazhdan-Lusztig polynomial
\cite[2.12]{A86},
which in turn is equal to
the ``leading" coefficient of
the relevant inverse
Kazhdan-Lusztig polynomial
\cite[8.2]{A86}.
One easily finds 
$\mu(B,A)=0$
by inspection
using
\cite{K87}.
The case when
$L_q(B)$ is a composition factor of
$\nabla_q(A)_3$ is analogous.
\end{pf}

\begin{rem}
The assumption that $\nabla_q(A)$ has simple head is satisfied for all $A \in \cA^{++}$ (Corollary 3.2). For type $B_2$ it is in fact satisfied for all $A \in \cA^+$
except for the alcoves numbered $9$ and $12$ below, see Appendix. Our direct computations below show that in the first case the Weyl module in question is rigid
whereas in 
the second case it is not. 
\end{rem}

\subsection{$B_2$-notation and -methods}
We let $\alpha_1$ denote the short simple root and $\alpha_2$ the long simple root. The corresponding reflections are denoted $s_1$ and $s_2$, respectively.
Then we set $S_l = \{s_0, s_1, s_2\}$ with $s_0$ denoting the reflection in the upper wall of $C$.

We enumerate the first few alcoves in the dominant chamber for $\mathrm{B}_2$ 
according to the strong linkage, see Fig.1.

\setlength{\unitlength}{1mm}
\begin{figure}
\begin{center}
\hspace*{5cm}
\begin{picture}(60,60)
\linethickness{1pt}
\put(0,0){\line(10,10){30}}
\put(0,0){\line(0,10){50}}
\put(0,10){\line(10,0){10}}
\put(10,10){\line(0,10){40}}
\put(10,10){\line(-10,10){10}}
\put(0,20){\line(10,0){20}}
\put(0,20){\line(10,10){30}}
\put(20,20){\line(0,10){30}}
\put(20,20){\line(-10,10){20}}
\put(0,30){\line(10,0){30}}
\put(30,30){\line(-10,10){20}}
\put(30,30){\line(0,10){20}}
\put(0,40){\line(10,0){30}}
\put(0,40){\line(10,10){10}}
\put(0,50){\line(10,0){30}}
\put(2,5){{\large 1}}
\put(2,12){{\large 2}}
\put(6,15){{\large 3}}
\put(6,22){{\large 5}}
\put(11,22){{\large 7}}
\put(11,15){{\Large 4}}
\put(2,26){{\large 6}}
\put(16,26){{\large 9}}
\put(20,26){{\large 11}}
\put(2,32){{\large 8}}
\put(15,32){{\large 12}}
\put(5.2,36){{\large 10}}
\end{picture}
\caption{Alcoves in type
$\mathrm{B}_2$}

\end{center}
\end{figure}
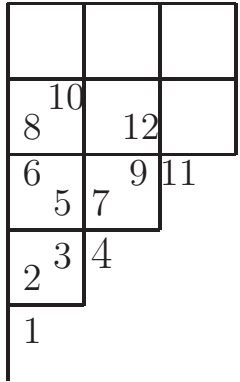

Using the above enumeration we write now $L_q(r)$, respectively $\Delta_q(r)$, respectively $T_q(r)$, short for the simple module $L_q(A_r)$, 
respectively the Weyl module $\Delta_q(A_r)$, respectively the indecomposable
tilting module $T_q(A_r)$, with $A_r$ being the alcove containing the number $r$.

If $M \in \cC_q$ we illustrate a filtration of $0 = M^{r+1} \subset M^r \subset \cdots M^1 \subset M^0 = M$ where $M^i/M^{i+1} =  
L_q(j^i_1) \oplus \cdots \oplus L_q(j^i_{n_i})$ by a diagram consisting of $r$ rows of boxes where the boxes in the $i$-th row contain
the numbers $j^i_1, \cdots, j^i_{n_i}$. The $r$-th row is at the bottom of the diagram. As a first step in determining such 
presentations of the $T_q(r)$'s we first consider their Weyl and dual Weyl filtrations. We picture those in a similar way by giving boxes containing the
appropriate modules but separated with dotted lines, meaning that the subquotient
admit further filtrations with those modules
as subquotients. 

In this way we will now picture 
the Loewy series of $\Delta_q(r)$ and $T_q(r), r = 1,, \cdots , 12$. It will turn out that all these modules are rigid except
$\Delta_q(12)$ and $T_q(12)$.
In fact, all $\Delta_q(r)$, $r>0$, are rigid except for
$r=12$.

Using the quantum Lusztig conjecture/theorem we can determine the composition factors of the Weyl modules and the corresponding tilting
conjecture/theorem gives us the Weyl filtrations of the indecomposable tilting modules. 

Our main tool for determining the Loewy structures 
of the Weyl modules will be the parity filtrations in Proposition 5.4 arising from the Loewy series of the corresponding baby Verma modules. 

In addition to this general method we shall often use translation arguments. In particular, we employ these to determine the socles of those Weyl modules where the 
Proposition in the Appendix does not apply. Observe that if for some $s \in S_l$ we have $A_r s < A_r$ then $\Delta_q(r) \subset \theta_s \Delta_q(r)$. Since 
$\theta_s L_q(i)$ = 0 for all $i$ with $A_i s < A_i$ we get from this 
\[ \text{ if $L_q(i) \subset \Delta_q(r)$ then $A_i s > A_i$ whenever $A_r s < A_r, s \in S_l$}.
\tag{1}
\]

\subsection{Alcove $1$}

$T_q(1)=\Delta_q(1)=L_q(1)$ is clearly rigid.
 
 \subsection{Alcove $2-6$}

The Weyl modules corresponding to these alcoves all have just two composition factors and the tilting modules just two Weyl factors (each having two composition factors).
This gives readily
 
$\Delta_q(2)=
\begin{tabular}{|c|}
\hline
$2$
\\
\hline
$1$
\\
\hline
\end{tabular}
$,
$\Delta_q(3)=
\begin{tabular}{|c|}
\hline
$3$
\\
\hline
$2$
\\
\hline
\end{tabular}
$,
$\Delta_q(4)=
\begin{tabular}{|c|}
\hline
$4$
\\
\hline
$3$
\\
\hline
\end{tabular}
$,
$\Delta_q(5)=
\begin{tabular}{|c|}
\hline
$5$
\\
\hline
$3$
\\
\hline
\end{tabular}
$, and 
$\Delta_q(6)=
\begin{tabular}{|c|}
\hline
$6$
\\
\hline
$5$
\\
\hline
\end{tabular}
$
are all rigid of Loewy length 2.

$T_q(2)=
\begin{tabular}{|c|}
\hline
$\Delta_q(1)$
\\
\hline
$\Delta_q(2)$
\\
\hline
\end{tabular}
=
\begin{tabular}{|c|}
\hline
$1$
\\
\hline
$2$
\\
\hline
$1$
\\
\hline
\end{tabular}
=
\begin{tabular}{|c|}
\hline
$\nabla_q(2)$
\\
\hline
$\nabla_q(1)$
\\
\hline
\end{tabular}
$, \;
$T_q(3)=
\begin{tabular}{|c|}
\hline
$\Delta_q(2)$
\\
\hline
$\Delta_q(3)$
\\
\hline
\end{tabular}
=
\begin{tabular}{|c|}
\hline
$2$
\\
\hline
$1\mid 3$
\\
\hline
$2$
\\
\hline
\end{tabular}
=
\begin{tabular}{|c|}
\hline
$\nabla_q(3)$
\\
\hline
$\nabla_q(2)$
\\
\hline
\end{tabular}
$, \;
$T_q(4)=
\begin{tabular}{|c|}
\hline
$\Delta_q(3)$
\\
\hline
$\Delta_q(4)$
\\
\hline
\end{tabular}
=
\begin{tabular}{|c|}
\hline
$3$
\\
\hline
$2\mid 4$
\\
\hline
$3$
\\
\hline
\end{tabular}
=
\begin{tabular}{|c|}
\hline
$\nabla_q(4)$
\\
\hline
$\nabla_q(3)$
\\
\hline
\end{tabular}
$, \;
$T_q(5)=
\begin{tabular}{|c|}
\hline
$\Delta_q(3)$
\\
\hline
$\Delta_q(5)$
\\
\hline
\end{tabular}
=
\begin{tabular}{|c|}
\hline
$3$
\\
\hline
$2\mid 5$
\\
\hline
$3$
\\
\hline
\end{tabular}
=
\begin{tabular}{|c|}
\hline
$\nabla_q(5)$
\\
\hline
$\nabla_q(3)$
\\
\hline
\end{tabular}$, and 
$T_q(6)=
\begin{tabular}{|c|}
\hline
$\Delta_q(5)$
\\
\hline
$\Delta_q(6)$
\\
\hline
\end{tabular}
=
\begin{tabular}{|c|}
\hline
$5$
\\
\hline
$3\mid 6$
\\
\hline
$5$
\\
\hline
\end{tabular}
=
\begin{tabular}{|c|}
\hline
$\nabla_q(6)$
\\
\hline
$\nabla_q(5)$
\\
\hline
\end{tabular}
$
are all rigid of Loewy length $3$.

\subsection{Alcove $7$}

The dual parity filtration $\Delta_q(7)^\bullet$
on $\Delta_q(7)$ coming from Proposition 5.4
reads
$\Delta_q(7) 
= 
\begin{tabular}{|c|}
\hline
$7$
\\
\hline
$5\mid4\mid2$
\\
\hline
$3$
\\
\hline
\end{tabular}$.
By 5.3(1) applied with $s = s_0$ and $s = s_1$ we see that 
$\soc\Delta_q(7)=L_q(3)$.
As the middle layer is semisimple by construction 
it follows that $\Delta_q(7)$ is rigid with the Loewy series given by the parity
$\Delta_q^\bullet$- filtration.

As $T_q(7)$ is a direct summand of
$\theta_{s_0} T_q(5)$, it follows from Corollary 3.15 that 
$ll T_q(7)\leq
ll T_q(5)+2=5$.
On the other hand, by Proposition 3.10 we have
$ ll T_q(7) \geq 2 ll \Delta_q(7) -1 = 5$, i.e., we have $ ll T_q(7) = 5$.

A $\Delta_q$- and a $\nabla_q$-filtration on
$T_q(7)$ may be expressed as 
$T_q(7)=
\begin{tabular}{|c|}
\hline
$\Delta_q(3)$
\\
\hline
$\Delta_q(4)\,\vdots\,\Delta_q(5)$
\\
\hline
$\Delta_q(7)$
\\
\hline
\end{tabular}
=
\begin{tabular}{|c|}
\hline
$\nabla_q(7)$
\\
\hline
$\nabla_q(4)\,\vdots\,\nabla_q(5)$
\\
\hline
$\nabla_q(3)$
\\
\hline
\end{tabular}$.

The first filtration shows that the possible simple factors of $\soc T_q(7)$ are $L_q(3)= \soc \Delta_q(7) = \soc\Delta_q(5)= 
\soc\Delta_q(4)$ and $ L_q(2) = \soc\Delta_q(3)$ whereas the second filtration gives the possibilities $L_q(r), \; r \in\{ 3, 4,5,7\}$. Hence
$L_q(3)$ is the only factor and it must occur with multiplicity $1$ because
\[
\Hom_{\cC_q}(L_q(3),T_q(7))
\subset
\Hom_{\cC_q}(L_q(3),\theta_{s_1}T_q(4))
\simeq
\Hom_{\cC_q}(\theta_{s_1}L_q(3),T_q(4)).
\]
is $1$-dimensional. In fact,  $ \theta_{s_1}L_q(3) = $
\begin{tabular}{|c|}
\hline
$3$
\\
\hline
$5$
\\
\hline
$3$
\\
\hline
\end{tabular}
with the projections onto $L_q(3)$ being the only homomorphisms into $T_q(4)$.

Next we determine $\soc_2 T_q(7)$. We have $\soc_2 \Delta_q(7) \subset \soc_2 T_q(7)$ and we
claim that in fact we have equality. Now any additional factor would belong to the socle of one of the other Weyl factors of $T_q(7)$, i.e., 
would be $L_q(2)$ or $L_q(3)$. The short exact sequence
$$ 0 \to \soc T_q(7) = L_q(3) \to T_q(7) \to T_q(7)/L_q(3) \to 0$$ first gives that $L_q(3)$ does
not occur in $\soc_2 T_q(7)$ (because 
\linebreak
$\Ext^1_{\cC_q}(L_q(3), L_q(3)) = 0$) and then that
$L_q(2)$ occurs at most once because $\Ext^1_{\cC_q}(L_q(2), L_q(3)) = \C$. 

Now the selfduality of $T_q(7)$ and the fact (observed above) that $ll T_q(7) = 5$ imply that the socle series is 
\[
T_q(7)=\begin{tabular}{|c|}
\hline
$3$
\\
\hline
$5\mid4\mid2$
\\
\hline
$3\mid7\mid3$
\\
\hline
$5\mid4\mid2$
\\
\hline
$3$
\\
\hline
\end{tabular}.
\]
This series is symmetric around the middle layer and hence the selfduality of $T_q(7)$ implies 
that it coincides with the radical series.

\subsection{Alcove $8$}
The same methods as used above give
that
$\Delta_q(8)$ is rigid of Loewy length 3,
with Loewy series
\[
\Delta_q(8)=
\begin{tabular}{|c|}
\hline
$8$
\\
\hline
$7\mid6$
\\
\hline
$5$
\\
\hline
\end{tabular}.
\]

Likewise we easily get that $ll T_q(8)$ has Loewy length $5$. It is rigid with Loewy structure 
\[
T_q(8)=\begin{tabular}{|c|}
\hline
$5$
\\
\hline
$6\ \mid7\mid\ 3$
\\
\hline
$8\mid5\mid4\mid5\mid2$
\\
\hline
$6\mid7\mid3$
\\
\hline
$5$
\\
\hline
\end{tabular}.
\]

\subsection{Alcove $9$}
The parity filtration on
$\Delta_q(9)$ is given by
$\Delta_q(9) 
= 
\begin{tabular}{|c|}
\hline
$9$
\\
\hline
$6\mid7\mid3\mid1$
\\
\hline
$5\ \mid\ 2$
\\
\hline
\end{tabular}$.
By 5.3(1) applied to $s_2$
we see that
$\soc\Delta_q(9)=L_q(5)\oplus
L_q(2)$.
As each layer is semisimple,
$ll \Delta_q(9)=3$.

Now suppose
$\rad^2\Delta_q(9)=L_q(5)$.
Then
$L_q(2)$ would lie in $\rad_1\Delta_q(9)$ and this contradicts the
parity vanishing 
$\Ext^1_{\cC_q}(L_q(2),L_q(9))=0$.
Thus,
$\rad^2\Delta_q(9)=L_q(5)\oplus
L_q(2)$.
It follows that both
the socle and the radical series of
$\Delta_q(9)$ equals its
parity filtration.

\begin{remark}
It turns out that this and $\Delta_q(12)$
are the only Weyl modules with non-simple socle
among all 
$\Delta_q(r)$'s, $r\geq1$, cf the Appendix.
\end{remark}

Now we look at $T_q(9)$. As $T_q(9)$ is a summand of $\theta_{s_2} T_q(7) =$
\linebreak
$\theta_{s_2} (\rad T_q(7)/\soc T_q(7))$ 
(from 5.7 we know the head and the socle of $T_q(7)$) it follows from Corollary 3.15 that $ll T_q(9) \leq 5$. By
Proposition 3.10 we must have equality. 

A $\Delta_q$-(resp.
$\nabla_q$-)
filtration of $T_q(9)$
reads
\[
T_q(9)=
\begin{tabular}{|c|}
\hline
$\Delta_q(2)$
\\
\hline
$\Delta_q(3)$
\\
\hline
$\Delta_q(5)$
\\
\hline
$\Delta_q(6)\,\vdots\,\Delta_q(7)$
\\
\hline
$\Delta_q(9)$
\\
\hline
\end{tabular}
=
\begin{tabular}{|c|}
\hline
$\nabla_q(9)$
\\
\hline
$\nabla_q(6)\,\vdots\,\nabla_q(7)$
\\
\hline
$\nabla_q(5)$
\\
\hline
$\nabla_q(3)$
\\
\hline
$\nabla_q(2)$
\\
\hline
\end{tabular}.
\]
The candidates for the factors of the socle of $T_q(9)$ are now $L_q(5), L_q(3)$, $L_q(2)$ (the intersection of the set of socles of the Weyl 
factors and the socles of the dual Weyl factors). However, we can erase $L_q(3)$ by the same method as in 5.3(1). This means
$\soc T_q(9)=\soc\Delta_q(9)=L_q(5)\oplus L_q(2)$.

Now Corollary 3.15 and the above easily give $ll \theta_{s_2} \Delta_q(7) = 4$. Hence we have 
\[
\theta_{s_2} \Delta_q(7)=\begin{tabular}{|c|}
\hline
$\Delta_q(7)$
\\
\hline
$\Delta_q(9)$
\\
\hline
\end{tabular}
=
\begin{tabular}{|c|}
\hline
$7$
\\
\hline
$9\mid5\mid4\mid2$
\\
\hline
$7\mid6\mid3\mid1\mid3$
\\
\hline
$5\quad\mid\quad2$
\\
\hline
\end{tabular}
\subset 
T_q(9).
\]
This means that
$\begin{tabular}{|c|}
\hline
$5\mid2$
\\
\hline
\end{tabular}
=\soc T_q(9)$,
$\begin{tabular}{|c|}
\hline
$7\mid6\mid3\mid1\mid3$
\\
\hline
\end{tabular}
\subset
\soc_2T_q(9)$,
\linebreak
$\begin{tabular}{|c|}
\hline
$9\mid5\mid4\mid2$
\\
\hline
\end{tabular}
\subset
\soc_3T_q(9)$,
and
$\begin{tabular}{|c|}
\hline
$7$
\\
\hline
\end{tabular}
\subset
\soc_4T_q(9)$.
On the other hand, looking at the dual Weyl filtration we see that 
\linebreak
$\Hom_{\cC_q}(L_q(3), T_q(9)/\soc T_q(9))$ is at most $2$-dimensional whereas all
\linebreak
$\Hom_{\cC_q}(L_q(1), T_q(9)/\soc T_q(9))$,  
$\Hom_{\cC_q}(L_q(7), T_q(9)/\soc T_q(9))$,
and 
\linebreak
$\Hom_{\cC_q}(L_q(6), T_q(9)/\soc T_q(9))$ are at most
$1$-dimensional. Also $L_q(2)$, $L_q(4)$, $L_q(5)$, and $L_q(6)$ cannot occur in $\soc_2 T_q(9)$ (they do not extend the socle).
It follows that
\[
\soc^2T_q(7)=
\begin{tabular}{|c|}
\hline
$7\mid6\mid3\mid1\mid3$
\\
\hline
$5\quad\mid\quad2$
\\
\hline
\end{tabular}.
\]

Recall that we have a homomorphism $\Delta_q(j) \to T_q(9)$ for each occurrence of $\nabla_q(j)$ in the $\nabla$-filtration of $T_q(9)$. Our results on the socle of $T_q(9)$
together with our findings in 5.5 show these homomorphisms in the case of $j =3$ and $j = 6$ both are injections. They then extend to embeddings $T_q(3)$
and $T_q(6)$
into
$T_q(9)$
(note that $T_q(3)$ and $T_q(6)$ have simple socles). By dualizing we get a surjection $T_q(9) \to T_q(3) \oplus T_q(6)$. It is clear that the submodule 
$\theta_{s_2} \nabla_q(7)$ is in the kernel of this surjection. By character considerations it is then equal to the kernel so that we have a short exact sequence
\[ 0 \to \theta_{s_2} \nabla_q(7) \to T_q(9) \to T_q(3) \oplus T_q(6) \to 0.
\]
It induces a commutative diagram of exact rows
\[
{\tiny
\xymatrix@C4mm{
0\ar[r]
&
\soc^3\theta_{s_2}\Delta_q(7)
\ar[d]
\ar[r]
&
\soc^3 T_q(9)
\ar[d]
\ar[r]
&
T_q(3)\oplus
T_q(6)
\ar[r]
&0
\\
0
\ar[r]
&
\soc^3\theta_{s_2}\Delta_q(7)/
\soc^2\theta_{s_2}\Delta_q(7)
\ar[r]
&
\soc^3 T_q(9)/\soc^2 T_q(9)
\ar[r]
&
\soc(T_q(3)\oplus
T_q(6)).
\ar@{^(->}!<0ex,3ex>;[u]
}
}
\]
As
$ll(T_q(9)/\soc^3T_q(9))=2$
while
$ll(T_q(3)\oplus T_q(6))=3$,
$\soc^3T_q(9)$ must surjects onto 
$\soc(T_q(3)\oplus
T_q(6))$,
and hence
\[
\soc_3T_q(9)=\begin{tabular}{|c|}
\hline
$9\mid5\mid4\mid2\mid2\mid5$
\\
\hline
\end{tabular}.
\]
It is now easy to complete the computations to get the following socle series for $T_q(9)$ 
\[
T_q(9)=
\begin{tabular}{|c|}
\hline
$5\quad\mid\quad2$
\\
\hline
$7\ \mid3\mid6\mid3\mid\ 1$
\\
\hline
$9\mid5\mid4\mid2\mid2\mid5$
\\
\hline
$7\ \mid6\mid3\mid1\mid\ 3$
\\
\hline
$5\quad\mid\quad2$
\\
\hline
\end{tabular}.
\]
The rigidity of
$T_q(9)$ follows from the self-duality and the symmetry of this socle series.

\subsection{Alcove $12$}

We claim that
the socle series of
$\Delta_q(12)$ is given by
\[
\begin{tabular}{|c|}
\hline
$12$
\\
\hline
$8\mid9\mid4\mid2$
\\
\hline
$7\ \mid6\mid\ 3$
\\
\hline
$5\ \mid\ 1$
\\
\hline
\end{tabular}
\]
and the radical series by
\[
\begin{tabular}{|c|}
\hline
$12$
\\
\hline
$8\mid9\mid4\mid2$
\\
\hline
$7\mid6\mid3\mid1$
\\
\hline
$5$
\\
\hline
\end{tabular}.
\]
The dual parity filtration
$\Delta_q^\bullet$ of
$\Delta_q(12)$ reads
\[
\Delta_q(12) = 
\begin{tabular}{|c|}
\hline
$12$
\\
\hline
$8\mid9\mid4\mid2$
\\
\hline
$7\mid6\mid3\mid1$
\\
\hline
$5$
\\
\hline
\end{tabular}.
\]
By Remark 5.2 we have $\rad^j\Delta_q(12)\subset \Delta_q(12)^j$ for all $j$.
As $L_q(5)$ lies in $\soc_4\nabla_q(12)$
by Lemma 3.20, it follows that
$ll\Delta_q(12)=4$
and that
$L_q(5)$ lies in $\rad_3\Delta_q(12)$,
i.e.
$
\rad_3\Delta_q(12)=L_q(5)$.
Now the structure of $\nabla_q(5)$ in 5.4  gives 
$\Ext^1_{\cC_q}(L_q(1),L_q(5))=0$ so that 
$L_q(1)\subset\soc\Delta_q(12)$.
On the other hand, 5.3(1) implies that $L_q(2), L_q(4), L_q(7)$, and $L_q(8)$ do not occur in $\soc\Delta_q(12)$.
Now the $\nabla_q$-filtration of $T_q(12)$ does not contain the factor $\nabla_q(3)$ and hence $\soc T_q(12)$ does not contain $L_q(3)$. Therefore
neither does $\soc \Delta_q(12)$. A direct computation of $\theta_{s_0} L_q(9)$ and $\theta_{s_0} L_q(6)$ reveals that $L_q(9)$ and $L_q(6)$ do
not map into $\Delta_q(12)$. So by checking all composition factors of $\Delta_q(12)$ we  have verified that 
\[\soc\Delta_q(12)=L_q(5)\oplus L_q(1).
\]

We have then 
$\soc^2\Delta_q(12)\supset
\begin{tabular}{|c|}
\hline
$7\mid6\mid3$
\\
\hline
$5\ \mid\ 1$
\\
\hline
\end{tabular}$.
We claim that this is an equality.
 
The structures of $\Delta_q(4)$ and $\Delta_q(8)$ which we have worked out in 5.5 and 5.7 show that 
$\Ext^1_{\cC_q}(L_q(8),L_q(j))=0=
\Ext^1_{\cC_q}(L_q(4),L_q(j))$, for $j \in \{1, 5\}$, see 5.2. This means that 
neither $L_q(8)$ nor $L_q(4)$ occur in $\soc_2\Delta_q(12)$.
Recall from Corollary 3.5 that there is a non-zero homomorphism $h$ from $\Delta_q(9)$ to $\Delta_q(12)$ and by Remark 3.6 this map becomes an
isomorphism when translated onto the $s_0$-wall. This means that the only possible composition factors of $\Ker h$ are $L_q(2)$ and $L_q(7)$. Therefore $\im h$
must have Loewy length $3$ so that $L_q(9)$ occurs in $\soc_3 \Delta_q(12)$. 

We claim that also $L_q(2)$ occurs in $\soc_3 \Delta_q(12)$. Otherwise, it must occur in $\soc_2 \Delta_q(12)$ and since $L_q(2)$ does not extend $L_q(5)$ this
means that (the only non-trivial extension of $L_q(1)$ and $L_q(2)$) $\Delta_q (2)$ must be contained in $\Delta_q(12)$. However, translating onto the $s_2$ wall 
this shows that we would have an injection
$\Delta_q(\mu) \to \Delta_q(\mu')$ where $\mu$ is a weight on the $s_2$-wall of $A_2$ and $\mu'$ similarly on the $s_2$-wall of $A_{12}$. This contradicts Corollary 3.2
and we have verified our claim.

Having now determined $\soc_2 \Delta_q(12)$ it only remains to observe that since $\hd \Delta_q(12) = L_q(12)$ we must have $\soc_3 \Delta_q(12) = \rad \Delta_q(12)$ so that 
the socle series is indeed the one stated.

Turning to the radical filtration, recall that we have already found that $\rad^3 \Delta_q(12) = L_q(5)$. By 5.1-2 of the remaining composition factors this 
leaves only the possibility 
$\rad_2\Delta_q(12)= $
\begin{tabular}{|c|}
\hline
$7\mid6\mid3\mid1$
\\
\hline
\end{tabular}
and 
$\rad_1\Delta_q(12)= $
\begin{tabular}{|c|}
\hline
$8\mid9\mid4\mid2$
\\
\hline
\end{tabular}
as stated.

\begin{remark}
The non-rigid Weyl module $\Delta_q(12)$ found here is in fact the only regular such Weyl module for type $B_2$. 
In the modular case this example was found by the first author in \cite{A87} but the full radical and socle series 
were not worked out.
\end{remark}

Consider now
$T_q(12)$. Its Weyl and dual Weyl filtrations are
\[
T_q(12)=
\begin{tabular}{|c|}
\hline
$\nabla_q(12)$
\\
\hline
$\nabla_q(9)\,\vdots\,\nabla_q(8)$
\\
\hline
$\nabla_q(7)\,\vdots\,\nabla_q(6)$
\\
\hline
$\nabla_q(5)$
\\
\hline
$\nabla_q(2)$
\\
\hline
$\nabla_q(1)$
\\
\hline
\end{tabular}
=
\begin{tabular}{|c|}
\hline
$\Delta_q(1)$
\\
\hline
$\Delta_q(2)$
\\
\hline
$\Delta_q(5)$
\\
\hline
$\Delta_q(7)\,\vdots\,\Delta_q(6)$
\\
\hline
$\Delta_q(9)\,\vdots\,\Delta_q(8)$
\\
\hline
$\Delta_q(12)$
\\
\hline
\end{tabular}.
\]
We 
shall limit ourselves to calculating enough of the Loewy structure to see that this is also a non-rigid module.

As for $T_q(9)$ we find that
$ll(T_q(12))=7$. Looking at the socles of the Weyl and dual Weyl modules occurring in the
filtrations of $T_q(12)$ limits the possible composition factors of $\soc T_q(12)$ to $L_q(1), L_q(2)$, and $L_q(5)$. Here the first and the last have
to occur as they do so in $\soc \Delta_q(12)$. On the other hand, 5.3(1) eliminates $L_q(2)$.
Thus,
$\soc T_q(12)=L_q(5)\oplus L_q(1)=\soc\Delta_q(12)$.

To determine $\soc_2 T_q(12)$ we note first that possible composition factors of 
$\soc_2T_q(12)/\soc_2\Delta_q(12)$ must be found among the socles of 
$\Delta_q(r)$, $r\in[1,9]\setminus\{3,4\}$.
In view of the vanishing extensions
with $\soc T_q(12)$ neither
$L_q(1)$ nor $L_q(5)$ can appear in
$\soc_2T_q(12)$. The same kind of reasoning shows that $L_q(3)$ and $L_q(2)$ both occur once in $\soc_2T_q(12)$. In fact, we get
\[
\Hom_{\cC_q}(L_q(3), T_q(12)/\soc T_q(12)) \subset \Ext^1_{\cC_q}(L_q(3), L_q(1) \oplus L_q(5)) = \C
\]
and an exact sequence
\[ \Hom_{\cC_q}(L_q(2), T_q(12)/\soc T_q(12)) \to \Ext^1_{\cC_q}(L_q(2), L_q(1) \oplus L_q(5)) \to \Ext^1_{\cC_q}(L_q(2), T_q(12)).
\]
Here the middle term is $\C$ and the last term is contained in 
\linebreak
$\Ext^1_{\cC_q}(L_q(2), \theta_{s_0} T_q(9)) \simeq \Ext^1_{\cC_q}(\theta_{s_0} L_q(2), T_q(9)) = 0$.
It follows that
$\soc_2T_q(12)=L_q(7)\oplus
L_q(6)\oplus
L_q(3)\oplus
L_q(2)
=
(\soc_2\Delta_q(12))\oplus
L_q(2)$,
where $L_q(2)$ comes from $\soc\Delta_q(9)$.
Thus
\[
\soc^3T_q(12)\supset
\begin{tabular}{|c|}
\hline
$8\mid9\mid4\mid2$
\\
\hline
$7\mid6\mid3\mid2$
\\
\hline
$5\ \mid\ 1$
\\
\hline
\end{tabular}.
\]
As $\Ext^1_{\cC_q}(L_q(5),L_q(12))=0$,
$\Delta_q(9)\subset T_q(12)/\Delta_q(12)$
has a summand
$L_q(5)$
of
$\soc\Delta_q(9)$
in
$\soc_3T_q(12)$
by parity.
To find other factors of
$\soc_3T_q(12)$
besides
$L_q(8)\oplus
L_q(9)\oplus
L_q(4)\oplus
L_q(2)\oplus
L_q(5)$, 
consider
$L_q(1)=\soc\Delta_q(2)$.
There is an exact sequence
\begin{multline*}
0\to\Hom_{\cC_q}(L_q(1),\soc^2T_q(12))\to
\Hom_{\cC_q}(L_q(1),T_q(12))
\to
\\
\Hom_{\cC_q}(L_q(1),T_q(12)/\soc^2 T_q(12))\to
\Ext^1_{\cC_q}(L_q(1),\soc^2 T_q(12))\to
\\
\Ext^1_{\cC_q}(L_q(1),T_q(12))
\end{multline*}
with
$\Ext^1_{\cC_q}(L_q(1),T_q(12))=0$
and 
$\Hom_{\cC_q}(L_q(1),\soc^2T_q(12))\simeq\C\simeq
\Hom_{\cC_q}(L_q(1),T_q(12))\simeq\Ext^1_{\cC_q}(L_q(1),\soc^2 T_q(12))$,
and hence
$[\soc_3T_q(12) : L_q(1)]=1$.

We note that the submodule 
$\Delta_q(12)$
extended by
$\Delta_q(9)$ in our $\Delta_q$-filtration of
$T_q(12)$
may be identified with
$\theta_0\Delta_q(9)$.
As
$4=ll(\Delta_q(12))\leq
ll(\theta_0\Delta_q(9))\leq
ll(\Delta_q(9))+2=5$
while
$\hd\theta_0\Delta_q(9)=L_q(9)$,
$ll(\theta_0\Delta_q(9))=5$.
Recalling that
the factor $L_q(1)$
in
$\soc_3T_q(12)$
comes from
$\Delta_q(2)$,
not from
$\Delta_q(9)$,
we see that component
$L_q(1)$ of $\Delta_q(9)$ must appear
in $\soc^4T_q(12)$.
It follows that
$[\soc^3T_q(12) : L_q(1)]=2$ while
\begin{align*}
1&\geq[T_q(12) : L_q(1)]-[\soc^4T_q(12) : L_q(1)]=
[T_q(12)/\soc^4T_q(12) : L_q(1)]
\\
&=
[\rad^4T_q(12) : L_q(1)]
\quad\text{by the self-duality of
$T_q(12)$},
\end{align*}
and hence
$\rad^4T_q(12)<\soc^3T_q(12)$,
verifying the nonrigidity of
$T_q(12)$.

\subsection{Weyl modules associated with higher alcoves}
The results in the appendix giving simple heads for certain dual Weyl modules (and thus simple socles of the corresponding Weyl modules)
show that for all the remaining alcoves $A \in \cA^+$ the socle of $\Delta_q(A)$ is simple. Combining with Corollary 5.6 we then get the rigidity of all these Weyl modules. Their
Loewy series therefore coincide with their parity filtration from Proposition 5.4, explicitly computable from the Loewy series of the corresponding baby 
Verma modules.

\subsection{Non-rigidity in the case of a singular weight}

Let $\mu$ be a weight in the $s_0$-wall and  denote by
$r'$ the image of $\mu$ in the the closure of
alcove
$r$.
We first show that
$\Delta_q(9')$ has
Loewy length 3
with socle series
\begin{tabular}{|c|}
\hline
$9'$
\\
\hline
$6'\mid3'$
\\
\hline
$5'\mid1'$
\\
\hline
\end{tabular}
while
the radical series is
\begin{tabular}{|c|}
\hline
$9'$
\\
\hline
$6'\mid3'\mid1'$
\\
\hline
$5'$
\\
\hline
\end{tabular}.

As $\Delta_q(9)$ is rigid with Loewy structure
\begin{tabular}{|c|}
\hline
$9$
\\
\hline
$7\mid6\mid3\mid1$
\\
\hline
$5\mid2$
\\
\hline
\end{tabular},
$\Delta_q(9') = T_\lambda^\mu \Delta_q(9)$
admits a filtration 
\begin{tabular}{|c|}
\hline
$9'$
\\
\hline
$6'\mid3'\mid1'$
\\
\hline
$5'$
\\
\hline
\end{tabular}
with semisimple subquotients.
If
$ll(\Delta_q(9'))=2$,
then
$L_q(6')\subset \soc\Delta_q(9')$,
and
hence
$L_q(6)\subset
T_\mu^\lambda
L_q(6')\subset
T_\mu^\lambda
\Delta_q(9')=
\begin{tabular}{|c|}
\hline
$\Delta_q(9)$
\\
\hline
$\Delta_q(12)$
\\
\hline
\end{tabular}$,
contradicting our findings in 5.7-8 giving
$\soc\Delta_q(9)=L_q(5)\oplus
L_q(2)$
and
$\soc\Delta_q(12)=L_q(5)\oplus
L_q(1)$.
Thus, the asserted radical series of
$\Delta_q(9')$ follows. 

On the other hand,
\begin{align*}
\Ext^1_{\cC_q}(L_q(1'),L_q(5'))
&\simeq
\Hom_{\cC_q}(L_q(1'),
\nabla_q(5')/L_q(5'))
\quad\text{as
$L_q(1')=\Delta_q(1')$}
\\
&\simeq
\Hom_{\cC_q}(L_q(1'),
L_q(3'))
=
0.
\end{align*}
Also,
as exhibited above,
$L_q(6')$ is not contained in $\soc\Delta_q(9')$,
and likewise we find that neither is 
$L_q(3')$.
Thus, the socle series of
$\Delta_q(9')$
follows as asserted.

Now recall from
\cite[5.2]{A00}
that
$T^\mu_\lambda
T_q(12)\simeq
T_q(9')^{\oplus2}$
and
$T^\mu_\lambda
T_q(7)\simeq
T_q(5')^{\oplus2}$.
It follows that
$T^\mu_\lambda
T_q(9)\simeq
T_q(9')\oplus
T_q(5')$
with
\[
T_q(9')=
\begin{tabular}{|c|}
\hline
$\Delta_q(1')$
\\
\hline
$\Delta_q(5')$
\\
\hline
$\Delta_q(6')$
\\
\hline
$\Delta_q(9')$
\\
\hline
\end{tabular}
\quad\text{and}\quad
T_q(5')=
\begin{tabular}{|c|}
\hline
$\Delta_q(3')$
\\
\hline
$\Delta_q(5')$
\\
\hline
\end{tabular}.
\]
By Proposition 3.10 we get 
$ll(T_q(9'))\geq 2ll(\Delta_q(9'))-1=5$.
On the other hand,
$T_q(9')$ is a direct summand of
$T^\mu_\lambda
T_q(9)$,
and hence
$ll(T_q(9'))\leq
ll(T_q(9))=5$,
forcing
$ll(T_q(9'))=5$.
The possible factors of
$\soc
T_q(9')$
are among
$L_q(9'),L_q(6'),L_q(5'),L_q(1')$
in view of the $\nabla_q$-filtration on
$T_q(9')$. Considering the socles of the factors in the $\Delta$-filtration of $T_q(9')$  we see that among these only 
$L_q(5')$ and $L_q(1')$ can occur. We conclude that 
$\soc T_q(9')=L_q(5')\oplus
L_q(1')=\soc\Delta_q(9')$.

To determine
$\soc_2T_q(9')$,
the possible factors besides those in
$\soc_2\Delta_q(9')$ are those in
$\soc\Delta_q(6'),
\soc\Delta_q(5'),
\soc\Delta_q(1')$.
As
$\Ext^1_{\cC_q}(L_q(5'),L_q(1'))=0$,
however,
neither 
$L_q(5')$ nor $L_q(1')$
can appear in
$\soc_2T_q(9')$.
Also,
\begin{multline*}
\Hom_{\cC_q}(L_q(3'),
T_q(9')/\soc T_q(9'))
\subset \Ext^1_{\cC_q}(L_q(3'),\soc T_q(9'))
\\
\simeq
\Ext^1_{\cC_q}(L_q(3'),L_q(5'))\oplus
\Ext^1_{\cC_q}(L_q(3'),L_q(1'))
\end{multline*}
with
$
\Ext^1_{\cC_q}(L_q(1'),L_q(3'))
=
\Ext^1_{\cC_q}(\Delta_q(1'),\nabla_q(3'))
=0$
while
\begin{align*}
\Ext^1_{\cC_q}(L_q(3'),L_q(5'))
&\simeq
\Hom_{\cC_q}(L_q(3'),\nabla_q(5')/L_q(5'))
\quad\text{as
$L_q(3')=\Delta_q(3')$}
\\
&\simeq
\Hom_{\cC_q}(L_q(3'),L_q(3'))
\simeq\C.
\end{align*}
Hence
$[\soc_2T_q(9') : L_q(3')]=1$ so that 
$\soc_2T_q(9')=
L_q(6')\oplus
L_q(3')=\soc_2\Delta_q(9')$,
and 
$\soc^2T_q(9')=\soc^2\Delta_q(9')$.
This implies
\begin{align*}
\rad^4T_q(9')&=\rad(\rad^3T_q(9'))
\subset \rad(\soc^2T_q(9'))
\quad\text{as
$ll(T_q(9'))=5$}
\\
&=
\rad(\soc^2\Delta_q(9'))
=\rad(\rad\Delta_q(9'))
\\
&=\rad^2\Delta_q(9')
\subset 
\soc\Delta_9(9')=\soc T_q(9')
\end{align*}
with the last containment being strict, see the above results on the Loewy structure on $\Delta_q(9')$.
This verifies the non-rigidity of
$T_q(9')$.

\begin{remark}
This example shows that $\Delta_q(A)$, respectively $T_q(A)$ may be rigid
without  the same being true for all
weights $\mu'$ in the closure of $A$ (take $A = A_9$ and $\mu' = 9'$). The same example also shows that $T_q(A)$ being rigid
is not enough to conclude the same for all summands of $\theta_s T_q(A)$ (take $A = A_9$ and note that the non-rigid module
(cf. 5.9) $T_q(12)$ is a summand of $\theta_{s_0} T_q(9)$). On the other hand, the example in 5.9 provides a non-rigid tilting
module, namely $T_q(12)$ which when we apply $\theta_{s_2}$ has only rigid summands. In fact, the Kazhdan-Lusztig theory gives $\theta_{s_2} T_q(12) = 
T_q(A_{(l-1)\rho}^+) \oplus T_q(9)$. The first summand is rigid by Theorem 3.27 and the last by 5.8. 
\end{remark}

\section{Non-rigid tilting modules for $SL_3(k)$}
In this section we take $G = SL_3(k)$. We shall give an explicit example of a non-rigid tilting module. The highest weight in question is only slightly above 
the upper bound in Theorem 4.7. So our example demonstrates that this bound cannot be relaxed and hence that the structure of modular tilting modules is
much more subtle than that of their quantized counterparts at complex roots of unity. At the same time our example relates to the recent work of Doty and Martin \cite{DM}. In fact,
we confirm their expectation that a certain module for $SL_3(k)$  with $k$ having characteristic $3$ is non-rigid, and we demonstrate that the same phenomena occur also
for all larger primes.

Our result is the following (using notation as in Section 4)

\begin{prop}
Let $G= SL_3(k)$ and assume $p$ is odd. Then $T(\frac{p(p-1)}{2} \rho)$ is non-rigid.
\end{prop}

\begin{remark} 
All composition factors $L(\lambda)$ of $T(\frac{p(p-1)}{2} \rho)$ except $L(\frac{p-1}{2} \rho)^{(1)}$ have $\lambda^1 \in \bar C$. Note that if $\alpha_1$ and $\alpha_2$ denote the two
simple roots then $\langle \frac{p-1}{2} \rho + \rho, (\alpha_1 + \alpha_2)^\vee \rangle = p+1$, i.e., $(\frac{p(p-1)}{2} \rho) ^1 = \frac{p-1}{2} \rho = 
s_0 \cdot (\frac{p-3}{2} \rho)$ lies in the second alcove $s_0\cdot C$ where
$s_0$ is the reflection in the upper wall of $C$. 
\end{remark}

\subsection{Proof for $p>3$}
Recall that for $G = SL_3(k)$ Conjecture 4.3 holds for all $p$ \cite[11.10]{RAG}. We denote by $Q(\lambda)$ the $G$-module such that $Q(\lambda)\vert_{G_1T} = Q_1(\lambda)$, $\lambda \in X_p$.
The isomorphism 4.7 can therefore in our case at hand be written 
$$T(\frac{p(p-1)}{2} \rho) \simeq Q((p-2) \rho) \otimes L(\frac{p-3}{2} \rho)^{(1)}.$$

In our case $h = 3$ so that Corollary 4.4 can be applied for $p > 5$. However, we claim that it also applies in the case $p =5$. In fact, in that case direct computations 
give that the relevant $\mu^1$'s are $0$, the two fundamental weights, and $\rho$. These all belong to the bottom alcove $C$ so that the proof of Corollary 4.4 still works, see
Remark 4.5. We get from Corollary 4.4 that $Q((p-2)\rho)$ is rigid of Loewy length $7$. Setting
$S^j = \soc_G^jQ((p-2)\rho)$ the assumptions in Lemma 4.6 are satisfied for all pairs $(S^7/S^j, L(\frac{p-3}{2} \rho))$ except when $ j = 3$. The middle layer $S^4/S^3$ contains
the
composition factor $L(p\rho)$ (with multiplicity $1$) and when tensored by $L(\frac{p-3}{2} \rho)^{(1)}$ the result is not semisimple. On the other hand, both $L(\rho)$ and 
$L(\frac{p-3}{2} \rho)$ are tilting modules (the highest weights of both come from $C$ since $p>3$). Therefore the tensor product $L(\rho) \otimes L(\frac{p-3}{2} \rho)$ 
is also tilting. The only non-simple summand of this tensor product is $T(\frac{p-1}{2} \rho)$ which equals the non-split extension 
$0 \rightarrow L(\frac{p-3}{2} \rho) \rightarrow T(\frac{p-1}{2} \rho) \rightarrow \nabla (\frac{p-1}{2} \rho ) \rightarrow 0$. 
It follows
that $T(\frac{p-1}{2} \rho)$ is uniserial of length $3$ with socle and head equal to $L(\frac{p-3}{2} \rho)$ and the middle factor equal to $L(\frac{p-1}{2} \rho)$.

Using the same approach as in the proof of Theorem 4.7 we get
$$ \soc_G^j T(\frac{p(p-1)}{2} \rho) = S^j \otimes L(\frac{p-3}{2} \rho)^{(1)} \text { for }  j = 1, 2, 3$$ 
and dually
$$ \rad_G^j T(\frac{p(p-1)}{2} \rho) = S^{7-j} \otimes L(\frac{p-3}{2} \rho)^{(1)} \text { for }  j = 1, 2, 3.$$
On the other hand, the above information shows that the layer $(S^4/S^3)\otimes L(\frac{p-3}{2} \rho)^{(1)}$ contains the composition factor $L(\frac{p(p-1)}{2} \rho)$ and
that this factor belongs to $ \rad_G^4T(\frac{p(p-1)}{2} \rho)$ but not to $\soc_G^4 T(\frac{p(p-1)}{2} \rho)$. If $T(\frac{p(p-1)}{2} \rho)$ were rigid then we would therefore have 
$ \rad_G^4 T(\frac{p(p-1)}{2} \rho) = \soc_G^j T(\frac{p(p-1)}{2} \rho)$ for some $j > 4$. However, $\soc_G^4 T(\frac{p(p-1)}{2} \rho)$ contains all the composition factors of the semisimple
part of $S^4/S^3 \otimes L(\frac{p-3}{2} \rho)^{(1)}$  and $ \rad_G^4 T(\frac{p(p-1)}{2} \rho)$ does not. This proves the non-rigidity of $T(\frac{p(p-1)}{2} \rho)$
in this case.

\subsection{The explicit socle series} 
In Fig. 2 we have numbered some of the alcoves in $X^+$ in such a way that the alcove containing $\frac{p(p-1)}{2} \rho$ has number $1$ and such that all other numbered
alcoves contains a weight for which the corresponding simple module is a composition factor in $T(\frac{p(p-1)}{2} \rho)$. Using the results from
\cite{AK} and the above observations,
we record below
the layers
$S^{j+1}/S^j \otimes L(\frac{p-3}{2} \rho)^{(1)}$ writing $m$ for the simple module with highest weight in the alcove numbered $m$. When $p = 5$ the alcoves numbered
$16$ and $19$ lie outside $X^+$ and those numbers should be ignored.

\[
\begin{tabular}{|c|}
\hline
$7$
\\
\hline
$13\quad\mid\quad8\quad\mid\quad11\quad\mid\quad18\quad\mid\quad6\quad\mid\quad15\quad\mid\quad17$
\\
\hline
$7^{\oplus2}\quad\mid\quad3\quad\mid\quad5\quad\mid\quad14\quad\mid\quad2\quad\mid\quad9\quad\mid\quad12$
\\
\hline
$13^{\oplus3}\mid8^{\oplus2}\mid11^{\oplus2}\mid18^{\oplus2}\mid6^{\oplus2}\mid15^{\oplus2}\mid17^{\oplus2}\mid4\mid20\mid10\mid16\mid19\mid T$
\\
\hline
$7^{\oplus2}\quad\mid\quad3\quad\mid\quad5\quad\mid\quad14\quad\mid\quad2\quad\mid\quad9\quad\mid\quad12$
\\
\hline
$13\quad\mid\quad8\quad\mid\quad11\quad\mid\quad18\quad\mid\quad6\quad\mid\quad15\quad\mid\quad17$
\\
\hline
$7$
\\
\hline
\end{tabular}
\]
with
$T=
T(\displaystyle\frac{p-1}{2}\rho)^{(1)}=
\begin{tabular}{|c|}
\hline
$13$
\\
\hline
$1$
\\
\hline
$13$
\\
\hline
\end{tabular}$.

\setlength{\unitlength}{1mm}
\begin{figure}
\begin{center}
\begin{picture}(150,100)
\linethickness{1pt}
\put(40,17.3){\line(100,0){60}}
\put(30,34.6){\line(100,0){80}}
\put(40,51.9){\line(100,0){60}}
\put(30,69.2){\line(100,0){80}}
\put(40,86.5){\line(100,0){60}}
\put(30,34.6){\line(100,173){40}}
\put(40,17.3){\line(100,173){40}}
\put(50,0){\line(100,173){50}}
\put(67,10){{\Large 20}}
\put(37,27){{\Large 16}}
\put(57,27){{\Large 17}}
\put(77,27){{\Large 18}}
\put(97,27){{\Large 19}}
\put(47,43){{\Large 11}}
\put(77,38){{\Large 14}}
\put(67,43){{\Large 13}}
\put(57,38){{\Large 12}}
\put(87,43){{\Large 15}}
\put(39,62){{\Large 4}}
\put(49,56){{\Large 5}}
\put(59,62){{\Large 6}}
\put(68,56){{\Large 7}}
\put(79,62){{\Large 8}}
\put(89,56){{\Large 9}}
\put(98,62){{\Large 10}}
\put(59,74){{\Large 2}}
\put(69,79){{\Large 1}}
\put(79,74){{\Large 3}}
\put(70,0){\line(100,173){40}}
\put(90,0){\line(100,173){20}}
\put(50,0){\line(-100,173){20}}
\put(70,0){\line(-100,173){40}}
\put(90,0){\line(-100,173){50}}
\put(100,17.3){\line(-100,173){40}}
\put(110,34.6){\line(-100,173){40}}
\end{picture}
\caption{
}

\end{center}
\end{figure}

\subsection{Proof in the case $p = 3$}

Finally, consider the case $p=3$, In this case we are dealing with $T(3 \rho) = Q(\rho)$. 
Since $T(3\rho)$ is selfdual, the middle layer 
$S = \soc_{G_1}^4 T(3\rho)/\soc^3_{G_1} T(3\rho)$ is the dual of $ \rad^3_{G_1} T(3\rho)/\rad^4_{G_1} T(3\rho) $ and by
the $G_1$-rigidity this equals $S$. In other words $S$ is also  selfdual. 
The socle series of
$Q(\rho)$ is given by
\[
\begin{tabular}{|c|}
\hline
$7$
\\
\hline
$6\quad\mid\quad 13\quad\mid\quad 8$
\\
\hline
$2\ \ \mid\  \ 7\ \  \mid\ \ 7\ \ \mid\ \ 3$
\\
\hline
$6\mid6\mid13 \mid8\mid8\mid T$
\\
\hline
$2\ \ \mid\  \ 7\ \  \mid\ \ 7\ \ \mid\ \ 3$
\\
\hline
$6\quad\mid\quad 13\quad\mid\quad 8$
\\
\hline
$7$
\\
\hline\end{tabular}
\quad\text{with
$T=\begin{tabular}{|c|}
\hline
$13$
\\
\hline
$1$
\\
\hline
$13$
\\
\hline
\end{tabular}$}.
\]

Consider now $Q= \{S/(L(\omega_1)^{\oplus2}\oplus
L(\omega_2)^{\oplus2})\}^{(-1)}$.
We claim that $Q$ is not semisimple. It it were then the composition factor $L(\rho)^{(1)}$ of $S$ would have to extend at least one 
factor of $\soc_G^3 T(\rho)/ \soc_G^2
T(\rho)$. The factors in question are 
$L(\rho + 3 \omega_1), L(\rho + 3 \omega_2)$ and $L(\rho)$. Now $\Ext_{G_1}^1(L(\rho), k) \simeq \Ext^1_{G_1}(k, L(\rho))$ can
be found by looking at the second socle layer in $T(3 \rho) = Q(\rho)$. Via Fig. 2 we get $\Ext^1_{G_1}(k, L(\rho))^{(-1)} = L(\omega_1) \oplus L(\omega_2)$. This
means that for any $\nu \in X^+$ we have  
\begin{align*}
\Ext_G^1(L(\rho)^{(1)}, L(\rho + 3\nu) 
&= \Hom_G(L(\rho), \Ext_{G_1}^1(k, L(\rho))^{(-1)} \otimes L(\nu)) 
\\
&= 
\Hom_G(L(\rho), (L(\omega_1) \oplus L(\omega_2)) \otimes L(\nu).
\end{align*}
This is certainly $0$ if $\nu = 0$. The other two relevant $\nu$'s are $\nu = \omega_1$ and $\nu =
\omega_2$. Easy calculations give $L(\omega_1) \otimes L(\omega_1) = L(2\omega) \oplus L(\omega_2)$ while $L(\omega_1) \otimes L(\omega_2) = T(\rho)$ and we
conclude that in all relevant cases the $\Ext^1$ vanishes.

We have shown that the middle layer of the $G_1$-Loewy series for $T(3 \rho)$ contains the non-semisimple $G$-summand $T(\rho)^{(1)}$. The same argument as used above
implies then that $T(3 \rho)$ is not rigid.

\section{Appendix}

Here we prove a general proposition which at least for low rank types is very
helpful in establishing that certain Weyl modules have simple socle, a result
that we need in Section 5. For convenience we formulate the proposition below only for
$G$ but both the statement and its proof carry over to  $U_q$ without any
change. A somewhat different argument for the same result can be found in 
Section 4.4 of \cite{A86a}.

We use the same notation as before and need a little more: We fix a Borel
subgroup $B$ in $G$ by requiring that the roots of $B$ are $-R^+$. If $\alpha
\in S$ then $P_\alpha$ denotes the minimal parabolic subgroup containing $B$
corresponding to $\alpha$. The induction functor from $B$-modules to
$G$-modules, respectively to $P_\alpha$-modules is written $H^0$, respectively
$H^0_\alpha$. The corresponding right derived functors are $H^i$ and
$H^i_\alpha$, respectively. If $\lambda \in X^+$ then $H^0(\lambda) =
\nabla(\lambda)$. Likewise we have for each $\lambda \in X$ satisfying $\langle
\lambda, \alpha^\vee \rangle \geq 0$ a non-zero $P_\alpha$-module $H^0_\alpha(\lambda)$ containing a unique
simple submodule which we denote $L_\alpha(\lambda)$.

\begin{prop}
Let $\lambda \in X^+$ satisfy  $\langle \lambda + \rho, \alpha^\vee \rangle < p$
and assume that $\lambda' = s_\alpha \cdot \lambda + p\alpha \in X^+$. Then
there is up to scalars a unique non-zero $G$-homomorphism $\nabla (\lambda') \to
\nabla(\lambda)$. Moreover, this homomorphism is surjective iff $H^2(s_\alpha
\cdot \lambda') = 0$.
\end{prop}

\begin{pf}
The $\lambda$ weight space $\nabla(\lambda')_\lambda$ is $1$-dimensional. Therefore there is
up to scalars at most one $G$-homomorphism from $\nabla(\lambda')$ to $\nabla(\lambda)$.
Easy $SL_2$-computations give the following two short exact sequences of
$P_\alpha$-modules
\[ 0 \to  L_\alpha(\lambda') \to H^0_\alpha(\lambda') \to H^0_\alpha(\lambda) \to 0
\tag{1} \]
and
\[ 0 \to H^1_\alpha(s_\alpha \cdot \lambda) \to 
H^1_\alpha(s_\alpha \cdot\lambda') \to L_\alpha(\lambda') \to 0.
\tag{2} \]
Moreover, $H^0_\alpha (\lambda) \simeq L_\alpha(\lambda)\simeq
H^1_\alpha(s_\alpha \cdot \lambda)$. 

Note that since $H^1_\alpha (\lambda) = 0$ we have $H^i(H^0_\alpha(\lambda))
\simeq
H^i(\lambda)$ and similarly for $\lambda'$. By Kempf's vanishing these vanish
for $i> 0$. Hence (1) gives the exact sequence
\[0 \to H^0(L_\alpha(\lambda')) \to \nabla(\lambda') \to \nabla(\lambda) \to 
H^1(L_\alpha(\lambda')) \to 0. \]

The isomorphism 
$H^0_\alpha (\lambda) \simeq H^1_\alpha(s_\alpha \cdot \lambda)$ combined with
Kempf's vanishing theorem give $H^i(H^1_\alpha(s_\alpha \cdot \lambda)) \simeq 
H^i(\lambda) = 0$ for $i>0$. Hence (2) implies $H^1(L_\alpha
(\lambda')) \simeq H^1(H^1_\alpha(s_\alpha \cdot \lambda'))$. But the vanishing
of $H^j_\alpha(s_\alpha \cdot \lambda')$ for $j \neq 1$ implies that 
$H^1(H^1_\alpha(s_\alpha \cdot \lambda')) \simeq H^2(s_\alpha \cdot \lambda')$.
Combining these facts we have proved the proposition.
\end{pf}

We now shift to the quantum case where we just add a subscript $q$ in the notation and replace $p$ by $l$. 
\begin{cor}
Let the notation and assumptions be as in the $q$-analogue of Proposition 7.1. Suppose in addition that $\lambda' \in l\rho + X^{+}$.
If $H^2_q(s_\alpha \cdot \lambda) = 0$ then $\nabla_q(\lambda)$ has simple head (equal to $\hd \nabla_q(\lambda')$).
\end{cor}

\begin{pf}
The proposition tells us that $\nabla_q(\lambda)$ is a quotient of $\nabla_q(\lambda')$. Hence the corollary
follows from Corollary 3.2.
\end{pf}

\begin{rem}
\begin{enumerate}
\item Suppose $R$ is of type $A_2$. We claim that in this case $\nabla_q(\lambda)$ has simple head for all $\lambda \in X^+$.
By corollary 3.2 it is enough to check this for $\lambda \in X^+\setminus (l\rho + X^{+})$. For such $\lambda$ we have $\lambda' \in X^+$ 
except if $\lambda$ belongs to one of the $4$ lowest alcoves in $X^+$ where the claim is easy to verify directly. For all remaining $\lambda$'s 
the relevant $H^2_q$ vanish (\cite {A81}) so that Corollary 7.2 applies.

\item Suppose $R$ is of type $B_2$. We claim that $\nabla_q(\lambda)$ has simple head for all $\lambda \in X^+\setminus Y$ where $Y$ is the union of alcoves 
$9$ and $12$, the wall between these and the wall between alcove $7$ and $9$ (we use the numbering from Figure 1). This is 
seen in exactly the same way as in (1). Note that if $\lambda \in Y$ then we have (with notation as above and with $\beta$ denoting
the long simple root) $H^2_q(s_\beta \cdot \lambda') \neq 0$, see \cite{A81}
\end{enumerate}
\end{rem}

\end{document}